\documentclass[12pt]{article}%
\usepackage{graphicx}
\usepackage{amsmath}
\usepackage{graphicx}
\usepackage{epstopdf}
\usepackage{geometry}
\usepackage{url}
\usepackage{amsthm}
\usepackage{amsxtra}
\usepackage{color}
\usepackage{lineno}
\usepackage{fancyhdr}%
\usepackage{amsfonts}%
\usepackage{amssymb}
\begin{document}

\begin{abstract}
We develop simple diffusion-advection models to estimate the average time it
takes fish to reach one of the boundaries of an enclosure and the population
distribution over time moving in the enclosure (such as a lake or slough). We
start with a combination of random walks and directed movement and then, from
these, proceeding to the associated Partial Differential Equations and their
solution. We also find the evolution of the population distribution and
communities composition over time moving in the enclosure. Although this model
was developed with fish movements in mind it has wide ranging applicability
scaling from the molecular to human and action from inert to deliberate.

\end{abstract}

\begin{titlepage}
\begin{center}
\vspace{2.5cm}  \title{Thinking Inside the Box: An Advection-Diffusion Model
of Animal Movement in an Enclosed Region}  \author{Stephen Tennenbaum*, John
Gatto and Joel Trexler \\ \\Department of Biology, Florida International
University \\ \\ \small{*set1@fiu.edu (corresponding author)}}     \maketitle%
\end{center}
\end{titlepage}

\renewcommand{\thefootnote}{\fnsymbol{footnote}}





\section{Introduction\medskip}

There is a wide range of movement phenomena occurring at scales from the
molecular to the human that share certain common features. There is the
dispersion of smoke from a chimney in the breeze, the spread of an effluent in
a stream, the transport of proteins in the cytoplasm of a cell, taxis of
protozoa in pond water responding to dissolved nutrients, movement of animals,
the spread of a tree species across a continent, these movements have both
directed components and random components. \ Both of these aspects can occur,
at one extreme, as completely passive or at the other, as an act of individual
volition. \ Nevertheless, if these two aspects of movement act at similar
scales in space and time then they can be described mathematically\ at any of
these scales using the similar principles. In this paper we derive an
advection-diffusion partial differential equation (PDE) from a mixed model
based on the hypothetical movements of fish in an enclosed body of water. This
model includes a random walk component corresponding to the fish searching its
immediate area and a directed movement component corresponding to the fish's
response to some environmental cue. The solution of this PDE allows us to
relate; (1) the parameters of movement to the average time to arrive at a
particular location from any starting point, (2) the probability of a fish
being at a particular location at a given time, \ and comparisons of the
movement characteristics of different species of fish (or different
individuals) with differing parameter values. This approach has been used in
studies fish movements in rivers, streams, and open water situations (see for
example \cite{Skalski00}, \cite{Sparrevohn02}, \cite{Faugeras07}, \cite{deKerckhove15}.) To
our knowledge this is the first time this method is applied to a completely
inclosed situation such as a lake, pond, marsh or slough.

We calculate first the average time to arrive at the designated goal starting
from any point in the enclosure. This will give us a general idea of how the
average individual performs given the particular aspects of its behavior. Next
we calculate the probability of being at any point in the enclosure at a
particular time. This will be useful for looking at the distribution of
individuals and any steady-state. \ We can also use this information to
determine from among different species the probabilities of first arrival at
the designated goal (or for that matter, any order of arrival). Finally, we
can determine the relative species composition of arrivals and distribution.

Our approach in both of the above cases will be to start with a random walk
alternated with directed movement. The derivations of diffusion type equations
from random walks can be found in many sources, but the one that we used here
is found in Random Walks in Biology \cite{Berg93}. \bigskip

\section{Average time to capture.}

\begin{figure}[th]
\centering\includegraphics[width=4.5in,height=3.3in]{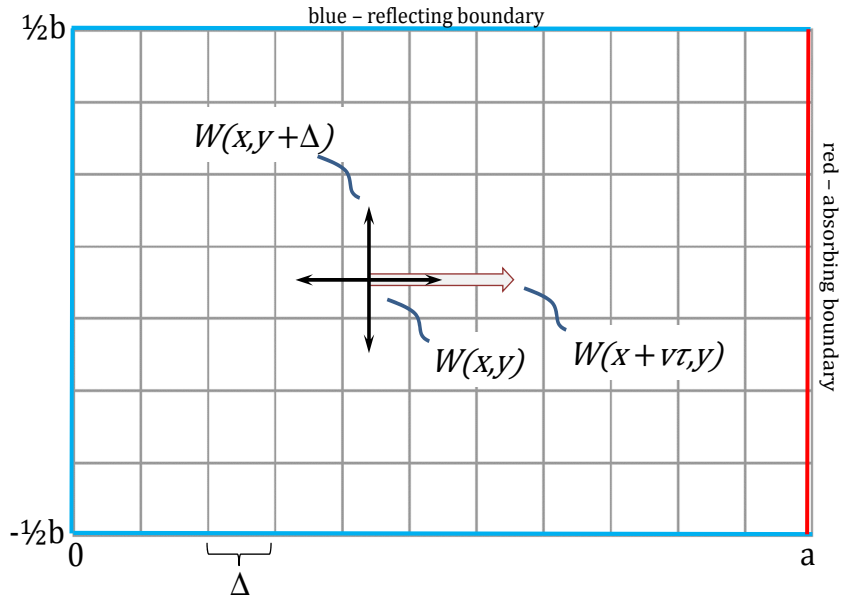}%
\caption{Mixed random walk and directed movement. Movement is directed a
fraction of the time $p$, and random movement a fraction of time $1-p = q$.
The average length of time to reach the absorbing boundary from a given point
is $W(x,y)$.}%
\label{Figure:1}%
\end{figure}

\subsection{Derivation of equations}

We start with the assumptions that movement will be either completely random
in two dimensions or completely deterministic in one dimension. The directed
movement will occur a fraction of time $p$ and the random movement will occur
a fraction of time $\left(  1-p\right)  .$ \ Random movement will involve
steps in only one of the four directions parallel to the boundaries of the
(rectangular) enclosure, each occurring with equal probability. Each of the
steps will move the individual a fixed distance $\Delta$ in a time $\tau$
(Figure \ref{Figure:1}). When the movement is directed, it will be exclusively
in the positive $x$ direction toward the designated goal. These steps shall
transport the individual a distance $v\tau$ toward the right, where $v$ is the
speed of directed movement and $\tau$ is our given time increment. \ Starting
at location $\left(  x,y\right)  ,$ we can write the average length of time to
get to the goal from that point as $W\left(  x,y\right)  .$ \ In the next time
increment $\tau$ the individual will have moved either randomly or directed.
If the movement was directed, it will now be at the location \ $\left(
x+v\tau,y\right)  $ and the average time it takes to get to goal from there
will be $W\left(  x+v\tau,y\right)  .$ \ Note that if there were only directed
movement then the time from $\left(  x,y\right)  $ would be $\tau$ units of
time more than from $\left(  x+v\tau,y\right)  ,$ which is closer to goal, so
that $W\left(  x,y\right)  =\tau+W\left(  x+v\tau,y\right)  .$ However,\ these
steps only occur with frequency $p,$ the rest of the time the individual is
moving in a random direction that will put it either closer $\left(
x+\Delta,y\right)  $ or further $\left(  x-\Delta,y\right)  $ or the same
distance but to one side or the other $\left(  x,y\pm\Delta\right)  $ in any
event the individual will be at these new locations $\tau$ units of time
later. \ So the average time from $\left(  x,y\right)  $ will be the average
of all the times from the new locations plus the time $\tau$ it took to get
there. \ Additionally, we assume that a fraction $s$ of the population does
not move at all in the interval, or alternatively, a given individual rests
$s$ fraction of the time. This is written explicitly as%

\begin{align}
W\left(  x,y\right)   &  =\tau+\left(  1-s\right)  \frac{1-p}{4}\left[
\smallskip W\left(  x+\Delta,y\right)  +W\left(  x-\Delta,y\right)  \right.
\nonumber\\
&  +\left.  W\left(  x,y+\Delta\right)  +W\left(  x,y-\Delta\right)
\smallskip\right]  +\left(  1-s\right)  pW\left(  x+v\tau,y\right)  +sW\left(
x,y\right)  . \label{1.01}%
\end{align}
\medskip Subtracting $W\left(  x,y\right)  $ from both sides, we
obtain\medskip%
\begin{align*}
0  &  =\tau+\left(  1-s\right)  \frac{1-p}{4}\left[  \smallskip\left(
W\left(  x+\Delta,y\right)  -W\left(  x,y\right)  \right)  -\left(  W\left(
x,y\right)  -W\left(  x-\Delta,y\right)  \right)  \right. \\
&  +\left.  \smallskip\left(  W\left(  x,y+\Delta\right)  -W\left(
x,y\right)  \right)  -\left(  W\left(  x,y\right)  -W\left(  x,y-\Delta
\right)  \right)  \right] \\
&  +\left(  1-s\right)  p\left(  W\left(  x+v\tau,y\right)  -W\left(
x,y\right)  \right)  .
\end{align*}
\medskip\ Then dividing through by $%
\frac12
\Delta,$ we get\medskip\
\begin{align}
0  &  =\frac{2\tau}{\left(  1-s\right)  \Delta}+\frac{1-p}{2}\left[  \left(
\frac{W\left(  x+\Delta,y\right)  -W\left(  x,y\right)  }{\Delta}\right)
-\left(  \frac{W\left(  x,y\right)  -W\left(  x-\Delta,y\right)  }{\Delta
}\right)  \right. \nonumber\\
&  +\left.  \left(  \frac{W\left(  x,y+\Delta\right)  -W\left(  x,y\right)
}{\Delta}\right)  -\left(  \frac{W\left(  x,y\right)  -W\left(  x,y-\Delta
\right)  }{\Delta}\right)  \right] \nonumber\\
&  +\frac{2pv\tau}{\Delta}\left(  \frac{W\left(  x+v\tau,y\right)  -W\left(
x,y\right)  }{v\tau}\right)  . \label{1.02}%
\end{align}
\medskip Noting that as $\Delta$ and $\tau$ approach $0,$ we get (by the
definition of the derivative)\medskip\
\begin{align*}
\lim\limits_{\Delta\rightarrow0}\frac{W\left(  x+\Delta,y\right)  -W\left(
x,y\right)  }{\Delta}  &  =\frac{\partial}{\partial x}W\left(  x+\Delta
,y\right) \\
\text{and \ }\lim\limits_{v\tau\rightarrow0}\frac{W\left(  x+v\tau,y\right)
-W\left(  x,y\right)  }{v\tau}  &  =\frac{\partial}{\partial x}W\left(
x+v\tau,y\right)  .
\end{align*}
\medskip\ Hence,\medskip%
\begin{align*}
0  &  =\frac{2\tau}{\left(  1-s\right)  \Delta}+\frac{1-p}{2}\left[  \left(
\frac{\partial}{\partial x}W\left(  x+\Delta,y\right)  -\frac{\partial
}{\partial x}W\left(  x,y\right)  \right)  \right. \\
&  +\left.  \left(  \frac{\partial}{\partial y}W\left(  x,y+\Delta\right)
-\frac{\partial}{\partial y}W\left(  x,y\right)  \right)  \right]
+\frac{2pv\tau}{\Delta}\frac{\partial}{\partial x}W\left(  x+\tau v,y\right)
.
\end{align*}
\medskip Dividing the terms in the bracket again by $\Delta,$ we obtain the
second derivative with\medskip\
\begin{align*}
&  \lim\limits_{\Delta\rightarrow0}\frac{\frac{\partial}{\partial x}W\left(
x+\Delta,y\right)  -\frac{\partial}{\partial x}W\left(  x,y\right)  }{\Delta
}=\frac{\partial^{2}}{\partial x^{2}}W\left(  x,y\right) \\
\text{and \ }  &  \lim\limits_{v\tau\rightarrow0}\frac{\partial}{\partial
x}W\left(  x+v\tau,y\right)  =\frac{\partial}{\partial x}W\left(  x,y\right)
.
\end{align*}
\medskip\ Thus, we continue as\medskip\
\[
0=\frac{2\tau}{\left(  1-s\right)  \Delta^{2}}+\frac{1-p}{2}\left(
\frac{\partial^{2}W}{\partial x^{2}}+\frac{\partial^{2}W}{\partial y^{2}%
}\right)  +\frac{2pv\tau}{\Delta^{2}}\frac{\partial W}{\partial x}.
\]
\medskip The diffusion coefficient is assumed to be a constant, for a 2d
random walk $D=\frac{\Delta^{2}}{4\tau}$ and $q=1-p,$ yielding\medskip%
\begin{align}
0  &  =\frac{1}{\left(  1-s\right)  2D}+\frac{q}{2}\left(  \frac{\partial
^{2}W}{\partial x^{2}}+\frac{\partial^{2}W}{\partial y^{2}}\right)  +\frac
{pv}{2D}\frac{\partial W}{\partial x}.\nonumber\\
0  &  =\frac{1}{\left(  1-s\right)  D}+q\left(  \frac{\partial^{2}W}{\partial
x^{2}}+\frac{\partial^{2}W}{\partial y^{2}}\right)  +\frac{pv}{D}%
\frac{\partial W}{\partial x}. \label{1.03a}%
\end{align}
\medskip Rearranging parameters and defining $\phi=\frac{pv}{qD}$ and
$\eta=\frac{1}{\left(  1-s\right)  qD},$we obtain\medskip%
\begin{equation}
0=\left(  \frac{\partial^{2}W}{\partial x^{2}}+\frac{\partial^{2}W}{\partial
y^{2}}\right)  +\phi\frac{\partial W}{\partial x}+\eta, \label{1.03b}%
\end{equation}
\medskip

The above equation does not lend itself to a solution by separation of
variables due to the presence of the constant $\eta$ (it is non-homogeneous),
so we will take as our ansatz a new function $U$ defined as,
\begin{equation}
U\left(  x,y\right)  =\phi W\left(  x,y\right)  +\eta x+k \label{1.04}%
\end{equation}
where $k$ is an arbitrary constant. \ We then have\medskip%
\[%
\begin{array}
[c]{ccc}%
\frac{\partial}{\partial x}U\left(  x,y\right)  =\phi\frac{\partial}{\partial
x}W\left(  x,y\right)  +\eta\smallskip & \text{and} & \frac{\partial^{2}%
}{\partial x^{2}}U\left(  x,y\right)  =\phi\frac{\partial^{2}}{\partial x^{2}%
}W\left(  x,y\right) \\
\frac{\partial}{\partial y}U\left(  x,y\right)  =\phi\frac{\partial}{\partial
y}W\left(  x,y\right)  & \text{and} & \frac{\partial^{2}}{\partial y^{2}%
}U\left(  x,y\right)  =\phi\frac{\partial^{2}}{\partial y^{2}}W\left(
x,y\right)
\end{array}
\]
\medskip Substituting the above into equation (\ref{1.03b}), we have\medskip\
\begin{equation}
0=\frac{1}{\phi}\left(  \frac{\partial^{2}U\left(  x,y\right)  }{\partial
x^{2}}+\frac{\partial^{2}U\left(  x,y\right)  }{\partial y^{2}}\right)
+\frac{\partial U\left(  x,y\right)  }{\partial x}. \label{1.05}%
\end{equation}
\medskip The next step is to try a separable solution for $U,$ that is,
$U\left(  x,y\right)  =X\left(  x\right)  Y\left(  y\right)  $ yields\medskip%

\[%
\begin{array}
[c]{ccc}%
\frac{\partial U\left(  x,y\right)  }{\partial x}=Y\left(  y\right)
\frac{\partial X\left(  x\right)  }{\partial x}\smallskip & \text{and} &
\frac{\partial^{2}U\left(  x,y\right)  }{\partial x^{2}}=Y\left(  y\right)
\frac{\partial^{2}X\left(  x\right)  }{\partial x^{2}}\\
\frac{\partial U\left(  x,y\right)  }{\partial y}=X\left(  x\right)
\frac{\partial\left(  y\right)  }{\partial y}Y & \text{and} & \frac
{\partial^{2}U\left(  x,y\right)  }{\partial y^{2}}=X\left(  x\right)
\frac{\partial^{2}Y\left(  y\right)  }{\partial y^{2}}%
\end{array}
\]
\medskip%
\[
0=\frac{1}{\phi}\left(  Y\frac{\partial^{2}X}{\partial x^{2}}+X\frac
{\partial^{2}Y}{\partial y^{2}}\right)  +Y\frac{\partial X}{\partial x}.
\]
\medskip Multiplying by $\phi$ and dividing by $XY,$ we get\medskip\
\begin{equation}
0=\left(  \frac{1}{X}\frac{\partial^{2}X}{\partial x^{2}}+\frac{1}{Y}%
\frac{\partial^{2}Y}{\partial y^{2}}\right)  +\phi\frac{1}{X}\frac{\partial
X}{\partial x}, \label{1.06}%
\end{equation}
\medskip and moving the $Y$ term to the left hand side gives\medskip%
\begin{equation}
-\frac{1}{Y}\frac{\partial^{2}Y}{\partial y^{2}}=\frac{1}{X}\left(
\frac{\partial^{2}X}{\partial x^{2}}+\phi\frac{\partial X}{\partial x}\right)
. \label{1.07}%
\end{equation}
\medskip Since each side of the equality in (\ref{1.07}) is independent of the
the other side, both sides must be equal to some constant (a standard PDE
approach)\medskip%
\begin{equation}
-\frac{1}{Y}\frac{\partial^{2}Y}{\partial y^{2}}=\lambda=\frac{1}{X}%
\frac{\partial^{2}X}{\partial x^{2}}+\phi\frac{1}{X}\frac{\partial X}{\partial
x}. \label{1.08}%
\end{equation}
Before attacking this equation, we discuss the boundary conditions. \ We have
three reflecting boundaries at $x=0,y=-b/2,$and $y=b/2.$ For these we write
$\frac{\partial}{\partial x}W\left(  0,y\right)  =0,$ and $\frac{\partial
}{\partial x}W\left(  x,\pm b/2\right)  =0$. \ At the goal, $x=a$ the boundary
condition is ``absorbing'', that is, the average time to reach the goal is
zero. \ So $W\left(  a,y\right)  =0.$ \ Going back to the definition of $U,$
\ we have\medskip%
\begin{align}
U\left(  a,y\right)   &  =\phi W\left(  a,y\right)  -\eta a+k\text{
\ \ \ \ from equation (\ref{1.04}) }\label{1.09}\\
X\left(  a\right)  Y\left(  y\right)   &  =-\eta a+k,\text{ \ \ \ since\ }%
U\left(  a,y\right)  =X\left(  a\right)  Y\left(  y\right)  \text{ and
}W\left(  a,y\right)  =0.\nonumber
\end{align}
\medskip

Solving for $Y,$ we obtain\medskip\
\begin{equation}
Y\left(  y\right)  =\frac{-\eta a+k}{X\left(  a\right)  }=C_{0},\text{ \ a
constant.} \label{1.10}%
\end{equation}
\medskip

Since the right hand side is composed entirely of constants and the function
$Y\left(  y\right)  $ is constant for all $x$ and $y,$ we have $\frac{\partial
Y}{\partial y}=0$ and $\frac{\partial^{2}Y}{\partial y^{2}}=0.$ \ Therefore,
$\lambda=0$ and we can now write\medskip\
\begin{equation}
\frac{\partial^{2}X}{\partial x^{2}}+\phi\frac{\partial X}{\partial x}=0.
\label{1.11}%
\end{equation}
\medskip

The associated characteristic equation is $r^{2}+\phi r=0,$ which has roots
$r_{1}=0$ and $r_{2}=-\phi$. This gives the general solution\medskip%
\begin{equation}
X\left(  x\right)  =C_{1}e^{r_{1}x}+C_{2}e^{r_{2}x}=C_{1}+C_{2}e^{-\phi x}
\label{1.12}%
\end{equation}
\medskip

Rearranging equation (\ref{1.04}) to solve for $W,$ we get\medskip%
\begin{align*}
W\left(  x,y\right)   &  =\phi^{-1}\left(  U\left(  x,y\right)  -\eta
x-k\right) \\
&  =\phi^{-1}\left(  X\left(  x\right)  Y\left(  y\right)  -\eta x-k\right)  .
\end{align*}
\medskip Substituting for $X$ from (\ref{1.12}) and $Y$ from (\ref{1.10}), we
get\medskip%
\begin{equation}
W\left(  x,y\right)  =\phi^{-1}\left(  \left(  C_{1}+C_{2}e^{-\phi x}\right)
C_{0}-\eta x-k\right)  . \label{1.13}%
\end{equation}
\medskip Taking the derivative with respect to $x,$ we have\medskip%
\begin{align}
\frac{\partial}{\partial x}W\left(  x,y\right)   &  =\frac{\partial}{\partial
x}\left(  \phi^{-1}\left(  \left(  C_{1}+C_{2}e^{-\phi x}\right)  C_{0}-\eta
x-k\right)  \right) \nonumber\\
&  =-C_{0}C_{2}e^{-\phi x}-\eta\phi^{-1}. \label{1.14}%
\end{align}
\medskip From our boundary condition at $x=0,$ we have $\frac{\partial
}{\partial x}W\left(  0,y\right)  =0.$ Plugging this into (\ref{1.14})
yields\medskip%
\begin{align}
\frac{\partial}{\partial x}W\left(  0,y\right)   &  =0=-\eta\phi^{-1}%
-C_{0}C_{2}e^{-\phi0}\nonumber\\
C_{0}C_{2}  &  =-\eta\phi^{-1}. \label{1.15}%
\end{align}
\medskip Equation (\ref{1.13}) then becomes%
\begin{align}
W\left(  x,y\right)   &  =\phi^{-1}\left(  \left(  C_{1}C_{0}+C_{2}%
C_{0}e^{-\phi x}\right)  -\eta x-k\right) \nonumber\\
&  =\phi^{-1}\left(  \left(  C_{1}C_{0}-\eta\phi^{-1}e^{-\phi x}\right)  -\eta
x-k\right) \nonumber\\
&  =\phi^{-1}\left(  \left(  C_{1}C_{0}-k\right)  -\eta\left(  \phi
^{-1}e^{-\phi x}+x\right)  \right)  . \label{1.16}%
\end{align}
\medskip At $x=a$ we have $W\left(  a,y\right)  =0,$ i.e.,\medskip\
\[
W\left(  a,y\right)  =0=\phi^{-1}\left(  \left(  C_{1}C_{0}-k\right)
-\eta\left(  \phi^{-1}e^{-\phi a}+a\right)  \right)
\]
\medskip Solving for the constants, we obtain\medskip%
\begin{equation}
C_{0}C_{1}-k=\eta\left(  \phi^{-1}e^{-\phi a}+a\right)  . \label{1.17}%
\end{equation}
\medskip Finally, substituting (\ref{1.17}) into the remaining constants in
(\ref{1.16}), we get\medskip%
\begin{align*}
W\left(  x,y\right)   &  =\phi^{-1}\left(  \left(  C_{0}C_{1}-k\right)
-\eta\left(  \phi^{-1}e^{-\phi x}+x\right)  \right) \\
&  =\phi^{-1}\left(  \eta\left(  \phi^{-1}e^{-\phi a}+a\right)  -\eta\left(
\phi^{-1}e^{-\phi x}+x\right)  \right) \\
&  =\eta\phi^{-1}\left(  a-x+\phi^{-1}e^{-\phi a}-\phi^{-1}e^{-\phi x}\right)
,\text{ \ or}%
\end{align*}
\medskip%
\begin{equation}
W\left(  x,y\right)  =\eta\phi^{-1}\left(  a-x\right)  -\eta\phi^{-2}\left(
e^{-\phi x}-e^{-\phi a}\right)  . \label{1.18}%
\end{equation}
\medskip\

Recalling that $\phi=\frac{pv}{qD}$ and $\eta=\frac{1}{\left(  1-s\right)
qD},$ we obtain\medskip%
\begin{equation}
W\left(  x,y\right)  =\frac{1}{\left(  1-s\right)  }\left(  \frac{\left(
a-x\right)  }{pv}-\frac{qD}{p^{2}v^{2}}\left(  e^{-\frac{pvx}{qD}}%
-e^{-\frac{pva}{qD}}\right)  \right)  . \label{1.19}%
\end{equation}
\medskip

Note that the solution is independent of the $y$ coordinate. \ On reflection
this makes perfect sense, since the back and side boundaries are reflecting,
the time to goal from any point on a line parallel to the goal should be the
same. This observation allows us to remark, in hindsight, that since $W\left(
x,y\right)  =W\left(  x\right)  ,$ then $\partial^{2}W\left(  x,y\right)
/\partial y^{2}=0$ with the result that the PDE could have been solved by
directly integrating twice and the exact same solution obtained.

Using equation (\ref{1.19}) we have the following special cases.

\begin{enumerate}
\item If the fish search and then move in a directed fashion very slowly then
$v\rightarrow0$ and $\lim\limits_{v\rightarrow0}W\left(  x,y\right)
=\frac{a^{2}-x^{2}}{\left(  1-s\right)  qD}.$

\item If they search all the time $q=1$ and don't rest $\left(  s=0\right)  ,$
then $W\left(  x,y\right)  =\frac{a^{2}-x^{2}}{2D}.$

\item If the fish move forward and then search very slowly for $q$ fraction of
the time, then $D\rightarrow0$ and $\lim\limits_{D\rightarrow0}W\left(
x,y\right)  =\frac{a-x}{\left(  1-s\right)  pv}.$

\item If they move forward without resting or searching, then $p=1,$ $s=0,$
and $W\left(  x,y\right)  =\frac{a-x}{v}$\bigskip
\end{enumerate}

This is not the only approach to this approach to this problem. We will see
later that the average time to arrive calculated here is slightly different
than the \textit{median }time to arrive calculated from the probability
distribution but has the advantage of being relatively simple to calculate. It
also allows us to characterize the system as dominated by random or directed
moment with a dimensionless index and can even provide a rough estimate of the
diffusion coefficient.

\begin{figure}[th]
\centering\includegraphics[width=4.5in,height=3.3in]{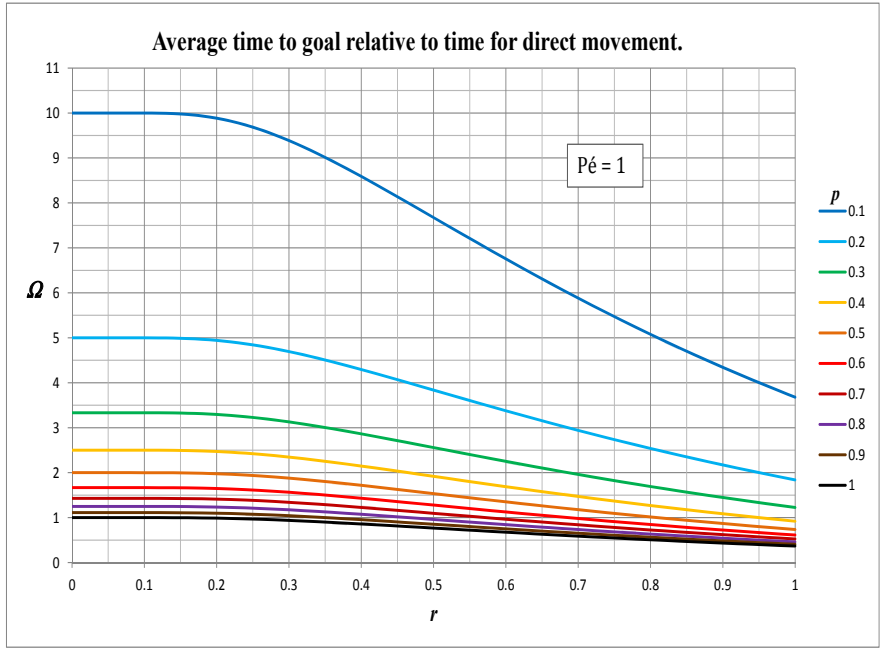}%
\caption{The time to reach the goal with random and constant directed movement
relative to the time to goal with just constant directed movement, $\Omega=
W(x,y)/((a-x)/v)$. Here, $p$ is the fraction of time spent in directed
movement, $r$ is the fraction of the distance from start to goal. Note that
for a fixed P\'{e}clet number the greater the percent time searching (i.e.,
diffusion), the longer the relative time to reach the goal on average. However
even if $p=1$ the time relative to advection exclusively is not $1$, since the
P\'{e}clet number includes divisor by $q$. Thus, in order for
${\mathrm{P\acute{e}}}=1$, the parameter $D \rightarrow\infty$.}%
\label{Figure:2}%
\end{figure}

For example, the relative importance of advection versus diffusion can be
expressed with the dimensionless quantity the \textit{P\'{e}clet number
}\cite{BCR12}. The P\'{e}clet number is defined as follows:%

\[
\mathrm{P\acute{e}}=\frac{\text{\textrm{advection transport rate}}%
}{\text{\textrm{diffusive transport rate}}}=\frac{v}{D/L}=\frac{vL}{D},
\]
\ where $L$ is the characteristic length of the system (for the above process
we will take it as the distance between the starting point and the goal on the
$x$-axis, i.e. $L=a-x$). \ For strictly physical processes this is fine,
however, in our formulation, we have included a parameter to indicate the
relative amount of time individuals do one thing or the other. \ Thus, we need
to include these modifications in the definition. We will do this by defining
an effective advection rate, $pv$ (directed speed), and effective diffusion
coefficient, $qD,$ so that
\[
\mathrm{P\acute{e}}=\frac{pv\left(  a-x\right)  }{qD}.
\]
\ Substituting this in the expression for the average time to goal allows us
to see the general behavior of the solution without having to manipulate
multiple variables or worrying about units of measurement. \ Rules of thumb
are for values of $\mathrm{P\acute{e}}<0.1$ diffusion dominates and for
$\mathrm{P\acute{e}}>10$ advection dominates. Thus, we can rewrite
(\ref{1.19}) as\medskip%
\begin{equation}
W\left(  x,y\right)  =\frac{a-x}{\left(  1-s\right)  pv}\left[  1-\frac
{1}{\mathrm{P\acute{e}}}\left(  e^{-\frac{x}{\left(  a-x\right)
}\mathrm{P\acute{e}}}-e^{-\frac{a}{\left(  a-x\right)  }\mathrm{P\acute{e}}%
}\right)  \right]  . \label{1.20}%
\end{equation}
\medskip Comparing $W$ to the time it takes to get to goal by directed
movement exclusively $\left(  \frac{a-x}{v}\right)  $ and letting
$r=1-\frac{x}{a}$ (the fraction of the distance to go to reach the goal) and
$s=0$ (no rest), we have\medskip%
\begin{equation}
\Omega=\frac{W\left(  x,y\right)  }{\left(  \frac{a-x}{v}\right)  }=\frac
{v}{\left(  \frac{a-x}{W\left(  x,y\right)  }\right)  }=\frac{1}{p}\left(
1-\frac{1}{\mathrm{P\acute{e}}}\left(  e^{-\frac{\left(  1-r\right)  }%
{r}\mathrm{P\acute{e}}}-e^{-\frac{1}{r}\mathrm{P\acute{e}}}\right)  \right)  .
\label{1.21}%
\end{equation}
\medskip This is the average time it takes to get to the goal relative to
traveling the entire distance at the advection rate. \ Alternatively, it is
the speed of advection exclusively relative to the average speed.A graph of
this for various values of $p$ is shown in Figure \ref{Figure:2}.\medskip

\medskip

The relationship in Equation (\ref{1.21}), relating the \textit{measurement}
of average time to goal relative to time required for direct travel and the
P\'{e}clet number, provides a possible means of estimating the diffusion
coefficient. For example, starting at the far end of the enclosure $\left(
r=1\right)  ,$ we have \medskip%
\begin{equation}
\Omega=\frac{1}{p}\left(  1-\frac{1}{\mathrm{P\acute{e}}}\left(
1-e^{-\mathrm{P\acute{e}}}\right)  \right)  . \label{1.22}%
\end{equation}
\medskip Graphing $\mathrm{P\acute{e}}$ as a function of $\Omega$ for various
values of $p,$ we have the following Figure \ref{Figure:2.1}.\medskip

\begin{figure}[th]
\centering\includegraphics[width=5in,height=2.76in]{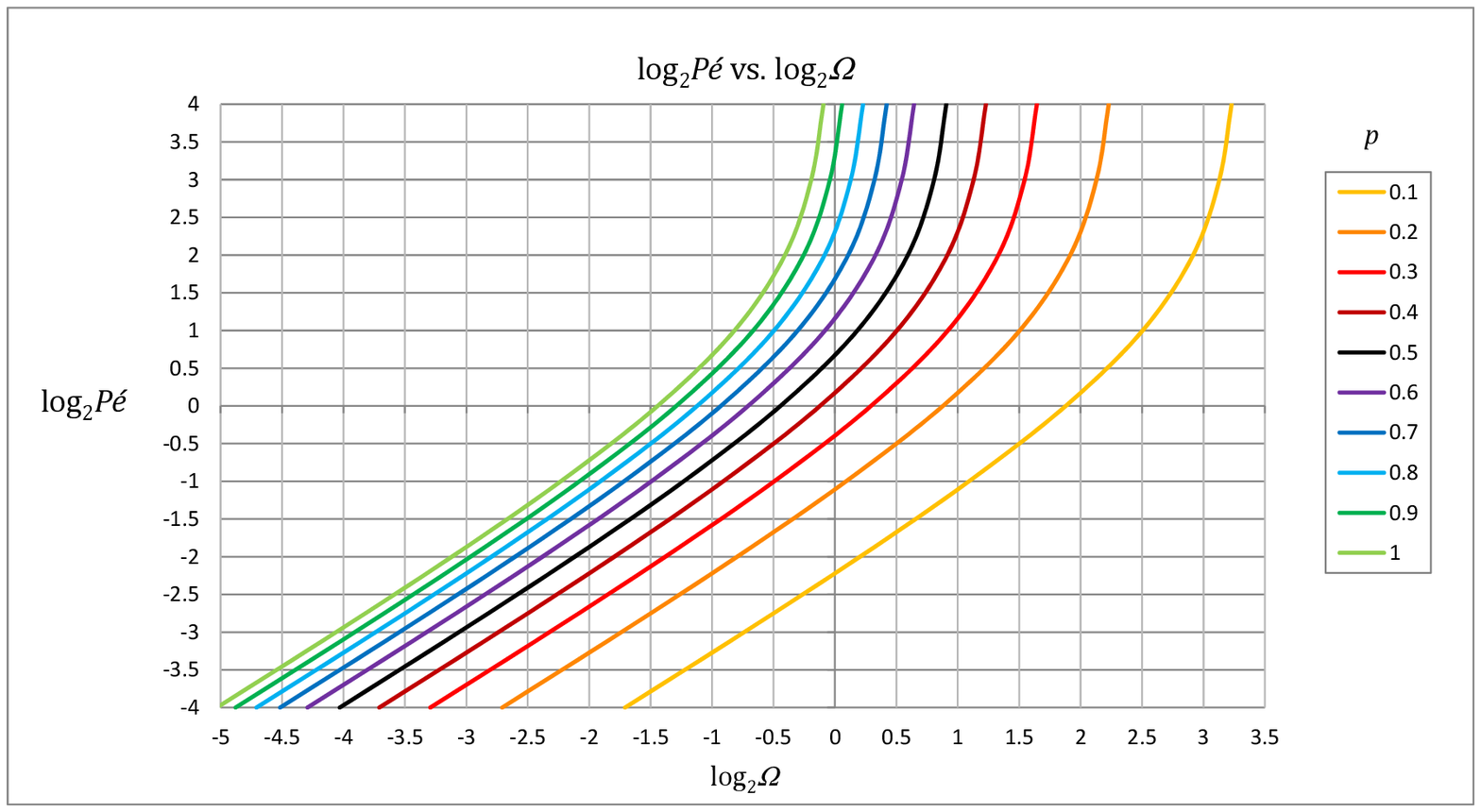}\caption{The
P\'{e}clet number as a function of time to goal relative to time to goal with
constant directed movement , plotted as the log base 2 of both variables
$log_{2}(P\acute{e})$ vs. $log_{2}(\Omega)$}%
\label{Figure:2.1}%
\end{figure}\medskip\bigskip

We now look at how the distribution of individuals actually evolves over time. \bigskip

\section{Progression of the probability distribution.}

\bigskip\begin{figure}[th]
\centering\includegraphics[width=5.5in,height=4.0in]{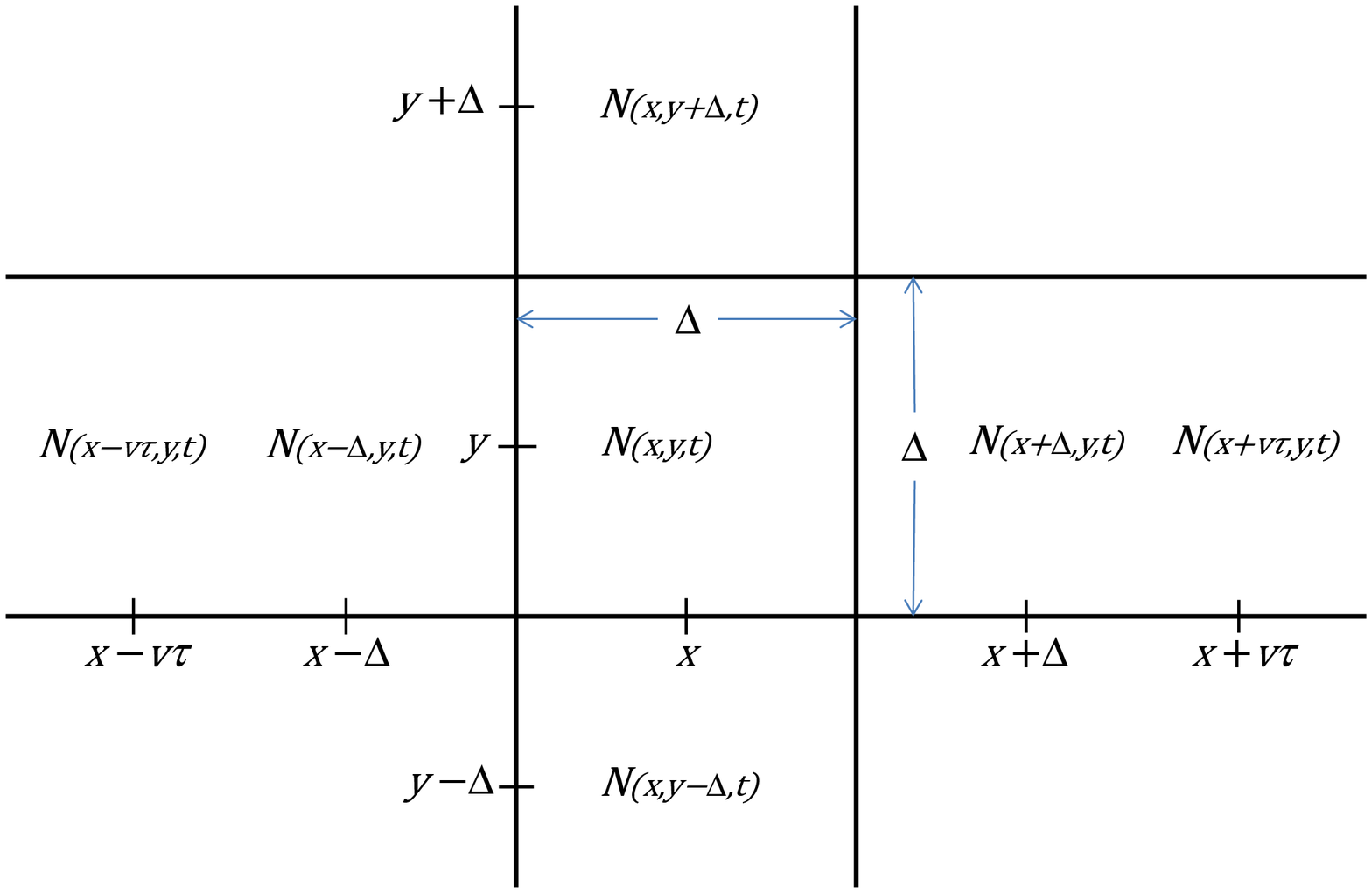}%
\caption{The number of individuals or particles relative to a point (x,y) at
time t.}%
\label{Figure:3}%
\end{figure}\medskip

In order to derive the 2-d probability distribution of individuals, we start
by looking at the change in numbers of individuals at a particular location.
The approach is similar to\ the above: we begin with the assumptions of a
semi-random walk in discrete time increments $\tau,$ moving in a directed
manner a fraction of the time $p$ and in a random manner a fraction of the
time $\left(  1-p\right)  =q$. \ On the right hand side we write the number of
individuals at location $\left(  x,y\right)  $ and at time $t+\tau.$ On the
left hand side is the number of individuals that have moved in the time
interval $\tau,$ from one location to another due to random or directed
movement. We assume that a fraction $s$ of all individuals everywhere will
\textit{not} have moved in the time interval $\tau$. With reference to Figure
3, we have\footnote{Alternatively we could formulate the problem in the
following way,
\begin{align*}
&  N\left(  x,y,t+\tau\right)  =\left(  1-s\right)  \left[  q^{\prime}N\left(
x+\Delta,y,t\right)  +p^{\prime}N\left(  x-\Delta,y,t\right)  \right. \\
&  \left.  +q^{\prime}N\left(  x,y+\Delta,t\right)  +q^{\prime}N\left(
x,y-\Delta,t\right)  \right]  +sN\left(  x,y,t\right)
\end{align*}
In this approach we have just a random walk but with a bias in the north
direction $p^{\prime}>q^{\prime}$ and $3q^{\prime}+p^{\prime}=1$. \ The
results of the derivation would be the same with $4q^{\prime}=q$ and
$p^{\prime}-q^{\prime}=p.$ The main difference is in the interpretation, in
the main text there is a random movement and a deterministic movement, in the
alternate derivation all movement is random but with a bias like dust in a
gentle breeze.}\medskip\
\begin{align}
&  N\left(  x,y,t+\tau\right)  =\left(  1-s\right)  \left[  \frac{{\small q}%
}{{\small 4}}N\left(  x+\Delta,y,t\right)  +\frac{{\small q}}{{\small 4}%
}N\left(  x-\Delta,y,t\right)  \right. \nonumber\\
&  \left.  +\frac{{\small q}}{{\small 4}}N\left(  x,y+\Delta,t\right)
+\frac{{\small q}}{{\small 4}}N\left(  x,y-\Delta,t\right)  +pN\left(
x-v\tau,y,t\right)  \right]  +sN\left(  x,y,t\right)  \label{2.01a}%
\end{align}
\medskip The fraction
${\frac14}$
in the above equation is because, for example, only a quarter of the number
moving from $\left(  x+\Delta,y\right)  $ go to $\left(  x,y\right)  ,$ the
other
${\frac34}$
go to one of the squares centered at $\left(  x+\Delta,y+\Delta\right)
,\left(  x+\Delta,y-\Delta\right)  ,$or $\left(  x+2\Delta,y\right)  .$ There
is no analogous (fraction) term on the directed movement, since the movement
is in one direction only. \ Subtracting $N\left(  x,y,t\right)  $ from both
sides, we have the change in numbers in the interval $\tau:$\medskip%
\begin{align}
&  N\left(  x,y,t+\tau\right)  -N\left(  x,y,t\right) \nonumber\\
&  =\left(  1-s\right)  \left[  \frac{{\small q}}{{\small 4}}\left(  N\left(
x+\Delta,y,t\right)  +N\left(  x-\Delta,y,t\right)  +N\left(  x,y+\Delta
,t\right)  +N\left(  x,y-\Delta,t\right)  \vspace{3pt}\right)  \right.
\nonumber\\
&  +\left.  pN\left(  x-v\tau,y,t\right)  \vspace{3pt}\right]  -\left(
1-s\right)  N\left(  x,y,t\right)  . \label{2.01b}%
\end{align}
\medskip Regrouping terms, we rewrite as\medskip%
\begin{align}
&  N\left(  x,y,t+\tau\right)  -N\left(  x,y,t\right) \nonumber\\
&  =\left(  1-s\right)  \frac{{\small q}}{{\small 4}}\left[  \left(  N\left(
x+\Delta,y,t\right)  -N\left(  x,y,t\right)  \right)  +\left(  N\left(
x-\Delta,y,t\right)  -N\left(  x,y,t\right)  \right)  \right. \nonumber\\
&  +\left.  \left(  N\left(  x,y+\Delta,t\right)  -N\left(  x,y,t\right)
\right)  +\left(  N\left(  x,y-\Delta,t\right)  -N\left(  x,y,t\right)
\right)  \right] \nonumber\\
&  +\left(  1-s\right)  p\left(  N\left(  x-v\tau,y,t\right)  -N\left(
x,y,t\right)  \right)  . \label{2.01c}%
\end{align}
\medskip Multiplying and dividing the random movement terms by the step
increment $\Delta,$ the directed movement term by the step increment $v\tau,$
and dividing both sides of the equation by the time increment $\tau,$ we
have\medskip\
\begin{align}
&  \frac{N\left(  x,y,t+\tau\right)  -N\left(  x,y,t\right)  }{\tau
}\nonumber\\
&  =\frac{\left(  1-s\right)  q\Delta}{4\tau}\left[  \left(  \frac{N\left(
x+\Delta,y,t\right)  -N\left(  x,y,t\right)  }{\Delta}\right)  -\left(
\frac{N\left(  x,y,t\right)  -N\left(  x-\Delta,y,t\right)  }{\Delta}\right)
\right. \nonumber\\
&  +\left.  \left(  \frac{N\left(  x,y+\Delta,t\right)  -N\left(
x,y,t\right)  }{\Delta}\right)  -\left(  \frac{N\left(  x,y,t\right)
-N\left(  x,y-\Delta,t\right)  }{\Delta}\right)  \right] \nonumber\\
&  -\frac{\left(  1-s\right)  pv\tau}{\tau}\left(  \frac{N\left(
x,y,t\right)  -N\left(  x-v\tau,y,t\right)  }{v\tau}\right)  . \label{2.01d}%
\end{align}
\medskip Taking the limit as $\tau,\Delta\rightarrow0$\medskip%
\begin{align}
&  \frac{\partial}{\partial t}N\left(  x,y,t+\tau\right)  =\frac{\left(
1-s\right)  q\Delta}{4\tau}\left[  \left(  \frac{\partial}{\partial x}N\left(
x+\Delta,y,t\right)  -\frac{\partial}{\partial x}N\left(  x,y,t\right)
\right)  \right. \nonumber\\
&  \left.  +\left(  \frac{\partial}{\partial y}N\left(  x,y+\Delta,t\right)
-\frac{\partial}{\partial y}N\left(  x,y,t\right)  \right)  \right]  -\left(
1-s\right)  pv\frac{\partial}{\partial x}N\left(  x,y,t\right)  .
\label{2.02a}%
\end{align}
\medskip\ \ Multiplying and dividing the random movement term by $\Delta$
again, we have\medskip%
\begin{align}
&  \frac{\partial}{\partial t}N\left(  x,y,t+\tau\right)  =\frac{\left(
1-s\right)  q\Delta^{2}}{4\tau}\left[  \left(  \frac{\frac{\partial}{\partial
x}N\left(  x+\Delta,y,t\right)  -\frac{\partial}{\partial x}N\left(
x,y,t\right)  }{\Delta}\right)  \right. \nonumber\\
&  \left.  +\left(  \frac{\frac{\partial}{\partial y}N\left(  x,y+\Delta
,t\right)  -\frac{\partial}{\partial y}N\left(  x,y,t\right)  }{\Delta
}\right)  \right]  -\left(  1-s\right)  pv\frac{\partial}{\partial x}N\left(
x,y,t\right)  . \label{2.02b}%
\end{align}
\medskip Taking the limit as $\tau\rightarrow0,\Delta\ \rightarrow0$
again\medskip%
\begin{equation}
\frac{\partial}{\partial t}N\left(  x,y,t\right)  =\frac{\left(  1-s\right)
q\Delta^{2}}{4\tau}\left(  \frac{\partial^{2}}{\partial x^{2}}N\left(
x,y,t\right)  +\frac{\partial^{2}}{\partial y^{2}}N\left(  x,y,t\right)
\right)  -\left(  1-s\right)  pv\frac{\partial}{\partial x}N\left(
x,y,t\right)  . \label{2.03}%
\end{equation}
\medskip

Recalling that the diffusion coefficient is $D=\frac{\Delta^{2}}{4\tau},$ we
get%
\begin{equation}
\frac{1}{\left(  1-s\right)  }\frac{\partial}{\partial t}N\left(
x,y,t\right)  =qD\left(  \frac{\partial^{2}}{\partial x^{2}}N\left(
x,y,t\right)  +\frac{\partial^{2}}{\partial y^{2}}N\left(  x,y,t\right)
\right)  -pv\frac{\partial}{\partial x}N\left(  x,y,t\right)  . \label{2.04}%
\end{equation}
\medskip

If we assume that the total population over the entire area in question is
constant and equal to $N_{T},$ then dividing through by this amount we have
the ``fraction''\ of the population at every point. \ The probability density
of finding an individual at a point is\ therefore $P\left(  x,y,t\right)
=N\left(  x,y,t\right)  /N_{T}$. Recall that\ in the average time to capture
calculation in the first section, the average time was scaled by the fraction
of time spent resting or staying put at a single location. \ Observe that the
same thing is occuring in (\ref{2.04}), rescaling time as $\left(  1-s\right)
t\rightarrow t$ (so that henceforth $t$ will represent the ``active'' time, or
the time for the active fraction of the population). \ If at some point in the
future we require the total time (say, for instance, we are measuring over the
course of a few days, or individuals are resting at regular intervals), then
$s$ becomes a significant fraction and can simply divide all $t$'s by $\left(
1-s\right)  $ to get the overall time.\ Equation (\ref{2.04}) simplifies
to\medskip%
\begin{equation}
\frac{\partial}{\partial t}P\left(  x,y,t\right)  =qD\left(  \frac
{\partial^{2}}{\partial x^{2}}P\left(  x,y,t\right)  +\frac{\partial^{2}%
}{\partial y^{2}}P\left(  x,y,t\right)  \right)  -pv\frac{\partial}{\partial
x}P\left(  x,y,t\right)  , \label{2.05}%
\end{equation}
where $t$ now represents movement durring active times.\bigskip

\subsection{Solution to diffusion-advection equation in a closed area.\medskip}

In order to solve this equation for $P,$ we will change our frame of reference
from the stationary coordinates $\left(  x,y,t\right)  $ to moving coordinates
$\left(  x^{\prime},y^{\prime},t^{\prime}\right)  $. The new reference frame
is moving at the speed $pv$ to the right (parallel to the stationary $x$-axis)
so that we have $x^{\prime}=x-\left(  x_{0}+pvt\right)  ,$ $y^{\prime}%
=y-y_{0},$ and $t^{\prime}=t.$ At $t=0,$ we have $\left(  x,y\right)  =\left(
x_{0},y_{0}\right)  $ and $\left(  x^{\prime},y^{\prime}\right)  =\left(
0,0\right)  .$ From these, we have the following derivatives\medskip%
\begin{align}
\frac{\partial x^{\prime}}{\partial t}  &  =-pv,\quad\frac{\partial x^{\prime
}}{\partial x}=1,\nonumber\\
\frac{\partial y^{\prime}}{\partial t}  &  =0,\quad\frac{\partial y^{\prime}%
}{\partial y}=1,\label{2.06}\\
\frac{\partial t^{\prime}}{\partial t}  &  =1,\quad\frac{\partial t^{\prime}%
}{\partial x}=0,\text{ \ \ and \ }\frac{\partial t^{\prime}}{\partial
y}=0.\nonumber
\end{align}
\medskip\ Using these and applying the chain rule we obtain the following
change of variables,\medskip\
\begin{align}
\frac{\partial}{\partial t}P\left(  x,y,t\right)   &  =\frac{\partial
}{\partial x^{\prime}}P\left(  x^{\prime},y^{\prime},t^{\prime}\right)
\frac{\partial x^{\prime}}{\partial t}+\frac{\partial}{\partial y^{\prime}%
}P\left(  x^{\prime},y^{\prime},t^{\prime}\right)  \frac{\partial y^{\prime}%
}{\partial t}+\frac{\partial}{\partial t^{\prime}}P\left(  x^{\prime
},y^{\prime},t^{\prime}\right)  \frac{\partial t^{\prime}}{\partial
t}\nonumber\\
&  =-pv\frac{\partial}{\partial x^{\prime}}P\left(  x^{\prime},y^{\prime
},t^{\prime}\right)  +\frac{\partial}{\partial t^{\prime}}P\left(  x^{\prime
},y^{\prime},t^{\prime}\right)  , \label{2.07a}%
\end{align}
\medskip%
\begin{align}
\frac{\partial}{\partial x}P\left(  x,y,t\right)   &  =\frac{\partial
}{\partial x^{\prime}}P\left(  x^{\prime},y^{\prime},t^{\prime}\right)
\frac{\partial x^{\prime}}{\partial x}+\frac{\partial}{\partial t^{\prime}%
}P\left(  x^{\prime},y^{\prime},t^{\prime}\right)  \frac{\partial t^{\prime}%
}{\partial x}=\frac{\partial}{\partial x^{\prime}}P\left(  x^{\prime
},y^{\prime},t^{\prime}\right) \label{2.07b}\\
\frac{\partial^{2}}{\partial x^{2}}P\left(  x,y,t\right)   &  =\frac{\partial
}{\partial x^{\prime}}\left(  \frac{\partial}{\partial x^{\prime}}P\left(
x^{\prime},y^{\prime},t^{\prime}\right)  \right)  \frac{\partial x^{\prime}%
}{\partial x}=\frac{\partial^{2}}{\left(  \partial x^{\prime}\right)  ^{2}%
}P\left(  x^{\prime},y^{\prime},t^{\prime}\right)  . \label{2.07c}%
\end{align}
\medskip Similarly,\medskip\
\begin{align}
\frac{\partial}{\partial y}P\left(  x,y,t\right)   &  =\frac{\partial
}{\partial y^{\prime}}P\left(  x^{\prime},y^{\prime},t^{\prime}\right)
\label{2.07d}\\
\frac{\partial^{2}}{\partial y^{2}}P\left(  x,y,t\right)   &  =\frac
{\partial^{2}}{\left(  \partial y^{\prime}\right)  ^{2}}P\left(  x^{\prime
},y^{\prime},t^{\prime}\right)  . \label{2.07e}%
\end{align}
\medskip Substituting the appropriate expressions from lines (\ref{2.06}) to
(\ref{2.07e}) into the derivatives in equation (\ref{2.05}), we have\medskip%
\begin{align*}
&  -pv\frac{\partial}{\partial x^{\prime}}P\left(  x^{\prime},y^{\prime
},t^{\prime}\right)  +\frac{\partial}{\partial t^{\prime}}P\left(  x^{\prime
},y^{\prime},t^{\prime}\right) \\
&  =qD\left(  \frac{\partial^{2}}{\left(  \partial x^{\prime}\right)  ^{2}%
}P\left(  x^{\prime},y^{\prime},t^{\prime}\right)  +\frac{\partial^{2}%
}{\left(  \partial y^{\prime}\right)  ^{2}}P\left(  x^{\prime},y^{\prime
},t^{\prime}\right)  \right)  -pv\frac{\partial}{\partial x^{\prime}}P\left(
x^{\prime},y^{\prime},t^{\prime}\right)
\end{align*}
\medskip Canceling the $\ -pv\frac{\partial}{\partial x^{\prime}}P\left(
x^{\prime},y^{\prime},t^{\prime}\right)  $ term on either side, we reduce
to\medskip\
\begin{equation}
\frac{\partial}{\partial t^{\prime}}P\left(  x^{\prime},y^{\prime},t^{\prime
}\right)  =qD\left(  \frac{\partial^{2}}{\left(  \partial x^{\prime}\right)
^{2}}P\left(  x^{\prime},y^{\prime},t^{\prime}\right)  +\frac{\partial^{2}%
}{\left(  \partial y^{\prime}\right)  ^{2}}P\left(  x^{\prime},y^{\prime
},t^{\prime}\right)  \right)  . \label{2.08}%
\end{equation}
\medskip This is just the form of an equation for diffusion in 2-d. \ We can
try separation of variables, where $P\left(  x^{\prime},y^{\prime},t^{\prime
}\right)  =X\left(  x^{\prime},t^{\prime}\right)  Y\left(  y^{\prime
},t^{\prime}\right)  ,$ to get\medskip%
\begin{align}
\frac{\partial}{\partial t^{\prime}}P\left(  x^{\prime},y^{\prime},t^{\prime
}\right)   &  =\frac{\partial}{\partial t^{\prime}}\left[  X\left(  x^{\prime
},t^{\prime}\right)  Y\left(  y^{\prime},t^{\prime}\right)  \right]
\nonumber\\
&  =Y\left(  y^{\prime},t^{\prime}\right)  \frac{\partial}{\partial t^{\prime
}}X\left(  x^{\prime},t^{\prime}\right)  +X\left(  x^{\prime},t^{\prime
}\right)  \frac{\partial}{\partial t^{\prime}}Y\left(  y^{\prime},t^{\prime
}\right)  , \label{2.09a}%
\end{align}%
\begin{align}
\frac{\partial^{2}}{\left(  \partial x^{\prime}\right)  ^{2}}P\left(
x^{\prime},y^{\prime},t^{\prime}\right)   &  =Y\left(  y^{\prime},t^{\prime
}\right)  \frac{\partial^{2}}{\left(  \partial x^{\prime}\right)  ^{2}%
}X\left(  x^{\prime},t^{\prime}\right)  ,\label{2.09b}\\
\text{and \ }\frac{\partial^{2}}{\left(  \partial y^{\prime}\right)  ^{2}%
}P\left(  x^{\prime},y^{\prime},t^{\prime}\right)   &  =X\left(  x^{\prime
},t^{\prime}\right)  \frac{\partial^{2}}{\left(  \partial y^{\prime}\right)
^{2}}Y\left(  y^{\prime},t^{\prime}\right)  . \label{2.09c}%
\end{align}
\medskip Substituting these in (\ref{2.08}) (and leaving off the function
arguments for brevity) gives us\medskip%
\begin{equation}
Y\frac{\partial X}{\partial t^{\prime}}+X\frac{\partial Y}{\partial t^{\prime
}}=qD\left(  Y\frac{\partial^{2}X}{\left(  \partial x^{\prime}\right)  ^{2}%
}+X\frac{\partial^{2}Y}{\left(  \partial y^{\prime}\right)  ^{2}}\right)  .
\label{2.10}%
\end{equation}
\medskip Putting all the $Y$ terms together and all the $X$ terms together we
have\medskip\
\begin{equation}
\frac{1}{Y}\left(  \frac{\partial Y}{\partial t^{\prime}}-qD\frac{\partial
^{2}Y}{\left(  \partial y^{\prime}\right)  ^{2}}\right)  +\frac{1}{X}\left(
\frac{\partial X}{\partial t^{\prime}}-qD\frac{\partial^{2}X}{\left(  \partial
x^{\prime}\right)  ^{2}}\right)  =0. \label{2.11}%
\end{equation}
\ \medskip Since both $X$ and $Y$ are positive functions and recalling that,
from the original construction of the problem, movement in the $x$ direction
is independent of movement in the $y$ direction, we determine that \medskip\
\begin{align}
\frac{\partial Y}{\partial t^{\prime}}-qD\frac{\partial^{2}Y}{\left(  \partial
y^{\prime}\right)  ^{2}}  &  =0\label{2.12a}\\
\text{and \ }\frac{\partial X}{\partial t^{\prime}}-qD\frac{\partial^{2}%
X}{\left(  \partial x^{\prime}\right)  ^{2}}  &  =0. \label{2.12b}%
\end{align}
\medskip We will start with a solution on an infinite domain\footnote{We could
also seek a solution directly on a finite domain for diffusion along the $y$
axis, however, the method used here is more succinct and provides a parallel
approach in both directions.}, i.e. \ $-\infty<x^{\prime}<\infty$ and
\ $-\infty<y^{\prime}<\infty$ \ and with the initial condition as a
instantaneous point source\medskip%
\begin{equation}
\left.
\begin{array}
[c]{c}%
Y\left(  y^{\prime},0\right)  =\delta\left(  y^{\prime}\right) \\
X\left(  x^{\prime},0\right)  =\delta\left(  x^{\prime}\right)
\end{array}
\right\}  \text{ \ Dirac delta functions.} \label{2.13}%
\end{equation}
\medskip Since the equations for $X$ and $Y$ have identical forms we will go
through the solution only for one of them. Applying the Fourier transform to
the $X$ equation, gives\medskip\
\begin{equation}
\mathcal{F}\left(  \frac{\partial}{\partial t^{\prime}}X\left(  x^{\prime
},t^{\prime}\right)  \right)  =qD\mathcal{F}\left(  \frac{\partial^{2}%
X}{\left(  \partial x^{\prime}\right)  ^{2}}\right)  . \label{2.14}%
\end{equation}
\medskip Which yields the ordinary differential equation (the circumflex over
the variable name indicates the variable is transformed from the time domain
to the frequency domain),\medskip\
\begin{equation}
\frac{d}{dt^{\prime}}\hat{X}\left(  \omega,t^{\prime}\right)  =-qD\omega
^{2}\hat{X}\left(  \omega,t^{\prime}\right)  . \label{2.15}%
\end{equation}
\medskip This has the solution\medskip%
\begin{equation}
\hat{X}\left(  \omega,t^{\prime}\right)  =\hat{X}\left(  \omega,0\right)
e^{-qDt^{\prime}\omega^{2}}. \label{2.16}%
\end{equation}
\medskip

We now need to return to the time domain. To do this we note that, \textit{in
general}, the Fourier transform of the convolution of two functions is the
product of their Fourier transforms. Therefore, the inverse Fourier transform
of a product of two functions in the frequency domain is the convolution of
the two functions in the time domain. For example,\medskip\
\begin{equation}
\mathcal{F}\left(  f\left(  x\right)  \ast g\left(  x\right)  \right)
=\mathcal{F}\left(  f\left(  x\right)  \right)  \mathcal{F}\left(  g\left(
x\right)  \right)  . \label{2.17}%
\end{equation}
\medskip Switching right and left-hand sides and taking the inverse transform
of both sides, we have\medskip\
\begin{equation}
\mathcal{F}^{-1}\left(  \mathcal{F}\left(  f\left(  x\right)  \right)
\mathcal{F}\left(  g\left(  x\right)  \right)  \right)  =\mathcal{F}%
^{-1}\left(  \mathcal{F}\left(  f\left(  x\right)  \ast g\left(  x\right)
\right)  \right)  =f\left(  x\right)  \ast g\left(  x\right)  , \label{2.18}%
\end{equation}
\medskip or, using the ``hat''notation for the functions in the frequency
domain ($\hat{f}\left(  \omega\right)  =\mathcal{F}\left(  f\left(  x\right)
\right)  $ and $\hat{g}\left(  \omega\right)  =\mathcal{F}\left(  g\left(
x\right)  \right)  $), we can rewrite the inverse transform of the product of
two transforms as\medskip\
\begin{equation}
\text{ }\mathcal{F}^{-1}\left(  \hat{f}\left(  \omega\right)  \hat{g}\left(
\omega\right)  \right)  =f\left(  x\right)  \ast g\left(  x\right)  .
\label{2.19}%
\end{equation}
\medskip

Employing the above device, we obtain the inverse transform in the following
steps.\medskip%
\begin{align}
X\left(  x^{\prime},t^{\prime}\right)   &  =\mathcal{F}^{-1}\left(  \hat
{X}\left(  \omega,t^{\prime}\right)  \right) \nonumber\\
&  =\mathcal{F}^{-1}\left(  \hat{X}\left(  \omega,0\right)  e^{-qDt^{\prime
}\omega^{2}}\right)  . \label{2.20}%
\end{align}
\medskip Note that\medskip%
\begin{align}
\mathcal{F}^{-1}\left(  \hat{X}\left(  \omega,0\right)  \right)   &  =X\left(
x^{\prime},0\right)  \Leftrightarrow\hat{X}\left(  \omega,0\right)
=\mathcal{F}\left(  X\left(  x^{\prime},0\right)  \right)  =\mathcal{F}\left(
\delta\left(  x^{\prime}\right)  \right) \label{2.21a}\\
\text{and \ }\mathcal{F}^{-1}\left(  e^{-qDt^{\prime}\omega^{2}}\right)   &
=\frac{1}{\sqrt{2qDt^{\prime}}}e^{-\frac{\left(  x^{\prime}\right)  ^{2}%
}{4qDt^{\prime}}}\Leftrightarrow e^{-\frac{qD}{2}t^{\prime}\omega^{2}%
}=\mathcal{F}\left(  \frac{1}{\sqrt{2qDt^{\prime}}}e^{-\frac{\left(
x^{\prime}\right)  ^{2}}{4qDt^{\prime}}}\right)  . \label{2.21b}%
\end{align}
\medskip Applying the convolution theorem for\ Fourier transforms and the
properties of the Dirac delta function, we get\medskip%
\begin{align}
X\left(  x^{\prime},t^{\prime}\right)   &  =\mathcal{F}^{-1}\left(  \hat
{X}\left(  \omega,0\right)  e^{-qDt^{\prime}\omega^{2}}\right) \nonumber\\
&  =\mathcal{F}^{-1}\left(  \mathcal{F}\left(  \delta\left(  x^{\prime
}\right)  \right)  \mathcal{F}\left(  \frac{1}{\sqrt{2qDt^{\prime}}}%
e^{-\frac{\left(  x^{\prime}\right)  ^{2}}{4qDt^{\prime}}}\right)  \right)
\nonumber\\
&  =\mathcal{F}^{-1}\left(  \mathcal{F}\left(  \delta\left(  x^{\prime
}\right)  \ast\frac{1}{\sqrt{2qDt^{\prime}}}e^{-\frac{\left(  x^{\prime
}\right)  ^{2}}{4qDt^{\prime}}}\right)  \right)  . \label{2.22}%
\end{align}
\medskip%
\begin{align}
X\left(  x^{\prime},t^{\prime}\right)   &  =\delta\left(  x^{\prime}\right)
\ast\frac{1}{\sqrt{2qDt^{\prime}}}e^{-\frac{\left(  x^{\prime}\right)  ^{2}%
}{4qDt^{\prime}}}\nonumber\\
&  =\frac{1}{\sqrt{4\pi qDt^{\prime}}}\int_{-\infty}^{\infty}\delta\left(
s\right)  e^{-\frac{\left(  x^{\prime}-s\right)  ^{2}}{4qDt^{\prime}}%
}ds\nonumber\\
&  =\frac{1}{\sqrt{4\pi qDt^{\prime}}}e^{-\frac{\left(  x^{\prime}\right)
^{2}}{4qDt^{\prime}}}. \label{2.23a}%
\end{align}
\medskip

With a similar result for the $Y$ function, we obtain\medskip%
\begin{equation}
Y\left(  y^{\prime},t^{\prime}\right)  =\frac{1}{\sqrt{4\pi qDt^{\prime}}%
}e^{-\frac{\left(  y^{\prime}\right)  ^{2}}{4qDt^{\prime}}}. \label{2.23b}%
\end{equation}
\medskip Switching back to the original reference frame we have\medskip%
\begin{align}
X\left(  x,t\right)   &  =\frac{1}{\sqrt{4\pi qDt}}e^{-\frac{\left(
x-x_{0}-pvt\right)  ^{2}}{4qDt}},\label{2.24a}\\
Y\left(  y,t\right)   &  =\frac{1}{\sqrt{4\pi qDt}}e^{-\frac{\left(
y-y_{0}\right)  ^{2}}{4qDt}}, \label{2.24b}%
\end{align}
\medskip

\ While we get solutions, this is not the end of the story. This result is
only valid on an unbounded plane but that is not the situation which we have
here. We have impenetrable boundaries at $y=\pm b/2$, $x=0$ and $x=a.$ The
above solution does not satisfy these conditions. To fulfill the desired
initial and boundary conditions we can add together multiple solutions to
create a new one by superposition given that the original differential
equation is linear. Three of the four boundaries $\left(  y=\pm b/2\text{ and
}x=0\right)  $ have no flux $\partial P\left(  0,y,t\right)  /\partial x=0,$
$\partial P\left(  x,b,t\right)  /\partial y=0,$ and $\partial P\left(
x,b,t\right)  /\partial x=0,$ they just reflect anything coming at them. The
fourth boundary, at $x=a,$ requires a bit more consideration, which we will
get to in a minute. For diffusion along the $y$-axis the fix is
straightforward, treat the two boundaries as if they were two mirrors facing
each other. We then take the part of each reflected image that lies between
$y=-b/2$ and $y=b/2$ and add it to the original distribution. \ If the
``$y$''\ distribution is centered at $y_{0}$ then the images are centered at
$nb+(-1)^{n}y_{0}$ for $n=\left\{  0,\pm1,\pm2,\pm3,\ldots\right\}  $. Or
letting $n_{\mathrm{even}}=2k$ for positive images and $n_{\mathrm{odd}}=2k+1$
for negative images, the $Y\left(  y,t\right)  $ factor is then\medskip%

\begin{equation}
\sum_{k=-\infty}^{\infty}Y_{k}\left(  y,t\right)  =\frac{1}{2\sqrt{\pi qDt}%
}\sum_{k=-\infty}^{\infty}\left(  e^{-\frac{\left(  y-y_{0}+2kb\right)  ^{2}%
}{4qDt}}+e^{-\frac{\left(  y+y_{0}+\left(  2k+1\right)  b\right)  ^{2}}{4qDt}%
}\right)  \text{ \ for }\ -b/2\leq y\leq b/2 \label{2.25a}%
\end{equation}
\medskip This spreads and becomes uniform distribution as time progresses, and
\ \medskip%
\begin{equation}
\lim_{t\rightarrow\infty}\sum_{k=-\infty}^{\infty}Y_{k}\left(  y,t\right)
=\frac{1}{b}. \label{2.25b}%
\end{equation}
\bigskip

\begin{figure}[th]
\centering\includegraphics[width=6in,height=2.2in]{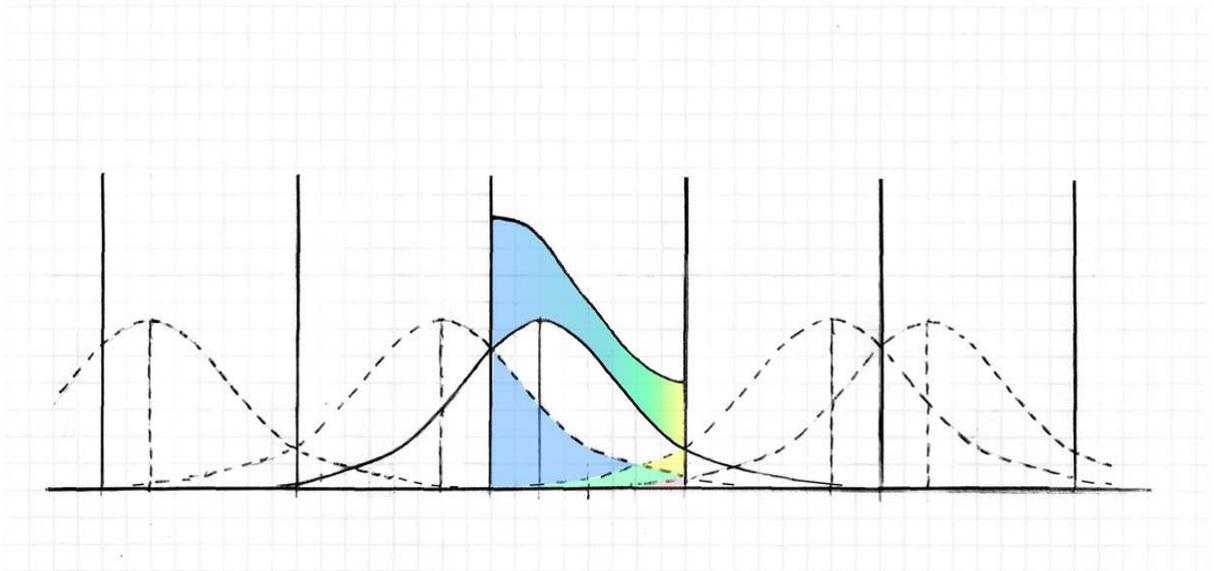}%
\caption{The primary distribution crosses the boundary at $y=-b/2$ and
$y=b/2$. These tails of the pdf can be exactly replaced within the boundaries
of the problem by adding a new ``reflected'' distributions centered at
$y=-b-y_{0}$ on the left and $y=b-y_{0}$ on the right. Each of these have
their own reflections at $y= \pm2b+y_{0}$ and so on ad infinitum. The
superposition of the imaged distributions (colored at bottom) in the practical
boundary is shown as blended colors over the primary distribution.}%
\label{Figure:4}%
\end{figure}\bigskip In practice only a few terms are needed since the tails
of the distributions coming from images that are far off from the original are
very small. Consequently the series converges quickly.

This approach does not work for movement parallel to the $x$-axis. \ If we
treat both of the boundaries at $x=0$ and $x=a$ as reflecting, the
distribution spreads with time as it does in diffusion only, but also the
entire distribution is translated to the right due to the advection term. This
looks alright until the center of the original distribution passes $x=a$. At
this point, since the distribution's reflection is traveling in the opposite
direction, more of the image distribution is between the boundaries traveling
to the left than the original still traveling to the right. Then as the first
image reaches the back wall it's image pokes it's head out and moves to the
right again. \ The entire affect is to make a superposed PDE that spreads as
it sloshes back and forth between the fore and aft walls!

What is to be done then? \ If we think about the situation of a chimney
emitting smoke into a very calm breeze the analogy is apt (see examples in
Chapter 2 sections 6 \& 7 of \cite{BCR12}): on the downwind side of the
chimney the breeze carries off the smoke, however, on the upwind side just
enough smoke diffuses in that direction to replace the amount carried back
downwind. If the amount of smoke and wind were perfectly unwavering, the plume
would quickly reach a steady state. \

Now, although this analogy is far from perfect, it does give us a starting
point for the analysis. At the upstream end, $x=0,$ we have a reflective
boundary $\frac{\partial}{\partial x}P\left(  0,y,t\right)  =0$ and at the
downstream end the barrier, instead of behaving as a boundary, acts like a
continuous source at $x=a,$ emitting particles into a steady flow, back into
the container by diffusion equal to the amount $\left(  Q\right)  $ to that
would have moved out if unimpeded. In other words, as the downstream end of
the original distribution moves mathematically past the boundary by a certain
amount, that amount actually diffuses back into the container. Again we start
with our advection-diffusion equation,
\begin{equation}
\frac{\partial}{\partial t}P\left(  x,y,t\right)  =qD\left(  \frac
{\partial^{2}}{\partial x^{2}}P\left(  x,y,t\right)  +\frac{\partial^{2}%
}{\partial y^{2}}P\left(  x,y,t\right)  \right)  -pv\frac{\partial}{\partial
x}P\left(  x,y,t\right)  . \label{3.01}%
\end{equation}
\medskip Since the flow is hitting a\ stationary wall, we apply separation of
variables without changing reference frame: $P\left(  x,y,t\right)
=\Psi\left(  x,t\right)  Y\left(  y,t\right)  $. We use $\Psi$ instead of $X$
to emphasize that this factor is a distribution due to the reflux parallel to
the $x$ axis and not the original flow. Thus we have\medskip%
\begin{align}
\frac{\partial}{\partial t}P\left(  x,y,t\right)   &  =\frac{\partial
}{\partial t}\left[  \Psi\left(  x,t\right)  Y\left(  y,t\right)  \right]
\nonumber\\
&  =Y\left(  y,t\right)  \frac{\partial}{\partial t}\Psi\left(  x,t\right)
+\Psi\left(  x,t\right)  \frac{\partial}{\partial t}Y\left(  y,t\right)
,\label{3.02a}\\
\frac{\partial}{\partial x}P\left(  x,y,t\right)   &  =Y\left(  y,t\right)
\frac{\partial}{\partial x}\Psi\left(  x,t\right)  ,\label{3.02b}\\
\frac{\partial^{2}}{\partial x^{2}}P\left(  x,y,t\right)   &  =Y\left(
y,t\right)  \frac{\partial^{2}}{\partial x^{2}}\Psi\left(  x,t\right)
,\label{3.02c}\\
\text{and \ }\frac{\partial^{2}}{\partial y^{2}}P\left(  x,y,t\right)   &
=\Psi\left(  x,t\right)  \frac{\partial^{2}}{\partial y^{2}}Y\left(
y,t\right)  . \label{3.02d}%
\end{align}
\medskip Substituting and regrouping like terms, we get\medskip%
\begin{align}
&  Y\left(  y,t\right)  \frac{\partial}{\partial t}\Psi\left(  x,t\right)
+\Psi\left(  x,t\right)  \frac{\partial}{\partial t}Y\left(  y,t\right)
\nonumber\\
&  =qD\left(  Y\left(  y,t\right)  \frac{\partial^{2}}{\partial x^{2}}%
\Psi\left(  x,t\right)  +\Psi\left(  x,t\right)  \frac{\partial^{2}}{\partial
y^{2}}Y\left(  y,t\right)  \right)  -pvY\left(  y,t\right)  \frac{\partial
}{\partial x}\Psi\left(  x,t\right)  ,\nonumber\\
0  &  =\frac{1}{\Psi\left(  x,t\right)  }\left(  \frac{\partial}{\partial
t}\Psi\left(  x,t\right)  -qD\frac{\partial^{2}}{\partial x^{2}}\Psi\left(
x,t\right)  +pv\frac{\partial}{\partial x}\Psi\left(  x,t\right)  \right)
\nonumber\\
&  +\frac{1}{Y\left(  y,t\right)  }\left(  \frac{\partial}{\partial t}Y\left(
y,t\right)  -qD\frac{\partial^{2}}{\partial y^{2}}Y\left(  y,t\right)
\right)  . \label{3.03}%
\end{align}
\medskip

The $Y$ terms are diffusion parallel to the $y$ axis, which we have solved for
above. The $\Psi$ term is advection and diffusion parallel to the $x$
axis\ but we now have a source term at the right hand boundary:\medskip%
\begin{equation}
\frac{\partial}{\partial t}\Psi\left(  x,t\right)  =qD\frac{\partial^{2}%
}{\partial x^{2}}\Psi\left(  x,t\right)  -pv\frac{\partial}{\partial x}%
\Psi\left(  x,t\right)  \label{3.04}%
\end{equation}
\medskip with initial and boundary conditions;\medskip%
\begin{align}
\Psi\left(  x,0\right)   &  =0\text{ \ \ \ initial conditions,}\label{3.05a}\\
\Psi\left(  -\infty,t\right)   &  =0,\text{\ \ left-hand side (momentarily
ignoring the boundary),}\label{3.05b}\\
\frac{\partial}{\partial t}\Psi\left(  a,t\right)   &  =\frac{d}{dt}Q\left(
t\right)  \text{ \ right-hand boundary (a variable source).} \label{3.05c}%
\end{align}
\medskip The Laplace transform in most useful in dealing with functions
defined on a semi-infinite domain. \ In preparation for applying a Laplace
transform we change coordinates (still fixed) from $x\in\left(  -\infty
,a\right)  $ to \ $x^{\prime}\in\left(  0,\infty\right)  ,$\medskip%
\begin{align}
\ x  &  =a-x^{\prime}\text{\ \ }\Leftrightarrow\text{ \ }x^{\prime
}=a-x\label{3.06a}\\
y  &  =y^{\prime}\text{ \ \ \ and \ \ \ }t=t^{\prime}\label{3.06b}\\
\frac{\partial t^{\prime}}{\partial t}  &  =1,\text{ \ \ }\frac{\partial
x^{\prime}}{\partial x}=-1,\text{ \ \ }\frac{\partial^{2}x^{\prime}}{\partial
x^{2}}=0, \label{3.06c}%
\end{align}
\medskip resulting in\medskip%
\begin{equation}
\frac{\partial}{\partial t^{\prime}}\Psi\left(  x^{\prime},t^{\prime}\right)
=\frac{qD}{2}\frac{\partial^{2}}{\left(  \partial x^{\prime}\right)  ^{2}}%
\Psi\left(  x^{\prime},t^{\prime}\right)  +pv\frac{\partial}{\partial
x^{\prime}}\Psi\left(  x^{\prime},t^{\prime}\right)  . \label{3.07}%
\end{equation}
\medskip Applying the Laplace transform (here we use a tilde to indicate the
transformed function, i.e., $\mathcal{L}\left(  \Psi\left(  x^{\prime
},t^{\prime}\right)  \right)  =\int_{0}^{\infty}\Psi\left(  x^{\prime
},t^{\prime}\right)  e^{-st}dt=\tilde{\Psi}\left(  x^{\prime},s\right)
$)\footnote{Here the ``$s$'' is the frequency variable of the Laplace
transform and not the fraction resting that we used in the random walk
example.},\medskip\
\begin{equation}
\mathcal{L}\left(  \frac{\partial}{\partial t^{\prime}}\Psi\left(  x^{\prime
},t^{\prime}\right)  \right)  =qD\mathcal{L}\left(  \frac{\partial^{2}%
}{\left(  \partial x^{\prime}\right)  ^{2}}\Psi\left(  x^{\prime},t^{\prime
}\right)  \right)  +pv\mathcal{L}\left(  \frac{\partial}{\partial x^{\prime}%
}\Psi\left(  x^{\prime},t^{\prime}\right)  \right)  , \label{3.08}%
\end{equation}
\medskip%
\begin{equation}
s\tilde{\Psi}\left(  x^{\prime},s\right)  -\Psi\left(  x^{\prime},0\right)
=qD\frac{d^{2}}{\left(  dx^{\prime}\right)  ^{2}}\tilde{\Psi}\left(
x^{\prime},s\right)  +pv\frac{d}{dx^{\prime}}\tilde{\Psi}\left(  x^{\prime
},s\right)  . \label{3.09}%
\end{equation}
\medskip Substituting the initial condition from (\ref{3.05a}), we get\medskip%
\begin{align}
s\tilde{\Psi}\left(  x^{\prime},s\right)   &  =qD\frac{d^{2}}{\left(
dx^{\prime}\right)  ^{2}}\tilde{\Psi}\left(  x^{\prime},s\right)  +pv\frac
{d}{dx^{\prime}}\tilde{\Psi}\left(  x^{\prime},s\right)  ,\label{3.10a}\\
0  &  =qD\frac{d^{2}}{\left(  dx^{\prime}\right)  ^{2}}\tilde{\Psi}\left(
x^{\prime},s\right)  +pv\frac{d}{dx^{\prime}}\tilde{\Psi}\left(  x^{\prime
},s\right)  -s\tilde{\Psi}\left(  x^{\prime},s\right)  . \label{3.10b}%
\end{align}
\medskip The transform of the boundary condition (\ref{3.05c}) gives\medskip%
\begin{align}
\Psi\left(  x^{\prime},0\right)   &  =0\text{ \ \ \ and \ \ \ \ }%
\frac{\partial}{\partial t^{\prime}}\Psi\left(  0,t^{\prime}\right)  =\frac
{d}{dt^{\prime}}Q\left(  t^{\prime}\right)  ,\nonumber\\
s\tilde{\Psi}\left(  0,s\right)  -\Psi\left(  0,0\right)   &  =s\tilde
{Q}\left(  s\right)  -Q\left(  0\right)  \Rightarrow\tilde{\Psi}\left(
0,s\right)  =\tilde{Q}\left(  s\right)  . \label{3.11}%
\end{align}
\medskip From (\ref{3.10b}) we have the characteristic equation\medskip%
\[
0=qDr^{2}+pvr-s,
\]
\medskip which has roots\medskip%
\begin{equation}
r_{1}=-\frac{pv}{2qD}\left(  1+\sqrt{1+\frac{4qDs}{p^{2}v^{2}}}\right)
,\text{ and\ \ }r_{2}=-\frac{pv}{2qD}\left(  1-\sqrt{1+\frac{4qDs}{p^{2}v^{2}%
}}\right)  . \label{3.12}%
\end{equation}
\medskip So the solution has the form\medskip%
\begin{equation}
\tilde{\Psi}\left(  x^{\prime},s\right)  =c_{1}e^{-\frac{pvx^{\prime}}%
{2qD}\left(  \sqrt{1+\frac{4qDs}{p^{2}v^{2}}}+1\right)  }+c_{2}e^{\frac
{pvx^{\prime}}{2qD}\left(  \sqrt{1+\frac{4qDs}{p^{2}v^{2}}}-1\right)  }.
\label{3.13}%
\end{equation}
\medskip Now as \ $s\rightarrow\infty,\ \tilde{\Psi}\left(  x^{\prime
},s\right)  <\infty\Longrightarrow c_{2}=0$ so that\medskip%
\begin{equation}
\tilde{\Psi}\left(  x^{\prime},s\right)  =c_{1}e^{-\frac{pvx^{\prime}}%
{2qD}\left(  \sqrt{1+\frac{4qDs}{p^{2}v^{2}}}+1\right)  }. \label{3.14}%
\end{equation}
\medskip However, at $x^{\prime}=0,\ \tilde{\Psi}\left(  0,s\right)  =c_{1},$
hence,\medskip\
\begin{equation}
\tilde{\Psi}\left(  x^{\prime},s\right)  =\tilde{\Psi}\left(  0,s\right)
e^{-\frac{pvx^{\prime}}{2qD}\left(  \sqrt{1+\frac{4qDs}{p^{2}v^{2}}}+1\right)
}. \label{3.15}%
\end{equation}
\medskip From (\ref{3.11}), we have%
\begin{equation}
\tilde{\Psi}\left(  x^{\prime},s\right)  =\tilde{Q}\left(  s\right)
e^{-\frac{pvx^{\prime}}{2qD}\left(  \sqrt{1+\frac{4qDs}{p^{2}v^{2}}}+1\right)
}. \label{3.16}%
\end{equation}
\medskip

We again need to return to the time domain, and as it was in the case of the
Fourier transform, the inverse Laplace transform of a product of two functions
in the frequency domain is the convolution of the two functions in the time
domain (or, equivalently, the convolution of the inverse Laplace transfoms of
the functions in the frequency domain) \footnote{The general approach,
starting with the Laplace transform of a convolution of two functions, is%
\[
\mathcal{L}\left(  f\left(  t\right)  \ast g\left(  t\right)  \right)
=\mathcal{L}\left(  f\left(  t\right)  \right)  \mathcal{L}\left(  g\left(
t\right)  \right)  .
\]
Switching the sides and applying the inverse transform, we have
\[
\mathcal{L}^{-1}\left(  \mathcal{L}\left(  f\left(  t\right)  \right)
\mathcal{L}\left(  g\left(  t\right)  \right)  \right)  =\mathcal{L}%
^{-1}\mathcal{L}\left(  f\left(  t\right)  \ast g\left(  t\right)  \right)
=f\left(  t\right)  \ast g\left(  t\right)  ,
\]
or using the ``tilda'' notation for the functions in the frequency domain we
have%
\[
\mathcal{L}^{-1}\left(  \tilde{f}\left(  s\right)  \tilde{g}\left(  s\right)
\right)  =f\left(  t\right)  \ast g\left(  t\right)  .
\]
}. \ Applying an inverse Laplace transform to both sides of (\ref{3.16}), we
have\medskip%
\begin{equation}
\Psi\left(  x^{\prime},t^{\prime}\right)  =\mathcal{L}^{-1}\left(  \tilde
{Q}\left(  s\right)  e^{-\frac{pvx^{\prime}}{2qD}\left(  \sqrt{1+\frac
{4qDs}{p^{2}v^{2}}}+1\right)  }\right)  . \label{3.17}%
\end{equation}
\medskip In (\ref{3.17}) the argument of the inverse Laplace transform is the
product of two Laplace transforms. \ Therefore the inverse transform is given
by\medskip\
\begin{equation}
\Psi\left(  x^{\prime},t^{\prime}\right)  =\mathcal{L}^{-1}\tilde{Q}\left(
s\right)  \ast\mathcal{L}^{-1}\left(  e^{-\frac{pvx^{\prime}}{2qD}\left(
\sqrt{1+\frac{4sqD}{p^{2}v^{2}}}+1\right)  }\right)  . \label{3.18}%
\end{equation}
\medskip Simplifying (\ref{3.18}) a bit, we arrive at\medskip\
\begin{equation}
\Psi\left(  x^{\prime},t^{\prime}\right)  =Q\left(  t^{\prime}\right)  \ast
e^{-\frac{pvx^{\prime}}{2qD}}\mathcal{L}^{-1}\left(  e^{-\frac{pvx^{\prime}%
}{2qD}\sqrt{1+\frac{4qD}{p^{2}v^{2}}s}}\right)  . \label{3.19}%
\end{equation}
\medskip In order to use Laplace transform tables to deal with the remaining
inverse transform, we need to simplify notation, which we can do by using the
following change of variables.\medskip\
\begin{align}
\xi &  =\frac{pvx^{\prime}}{2qD},\text{ \ \ and \ \ }\gamma=1+\frac{4qD}%
{p^{2}v^{2}}s,\label{3.20a}\\
s  &  =\frac{p^{2}v^{2}}{4qD}\left(  \gamma-1\right)  \Rightarrow
ds=\frac{p^{2}v^{2}}{4qD}d\gamma. \label{3.20b}%
\end{align}
\medskip We write out the inverse transform explicitly to sort out the
re-scaling necessary when changing variables:\medskip\
\begin{equation}
\mathcal{L}^{-1}\left(  e^{-\frac{pvx^{\prime}}{2qD}\sqrt{1+\frac{4qD}%
{p^{2}v^{2}}s}}\right)  =\frac{1}{2\pi i}\int_{c-i\infty}^{c+i\infty}%
e^{-\frac{pvx^{\prime}}{2qD}\sqrt{1+\frac{4qD}{p^{2}v^{2}}s}}e^{st^{\prime}%
}ds. \label{3.21}%
\end{equation}
\medskip Now substituting $\ x^{\prime}\rightarrow\xi$ and $s\rightarrow
\gamma$ terms, we get\medskip%
\begin{align}
\mathcal{L}^{-1}\left(  e^{-\frac{pvx^{\prime}}{2qD}\sqrt{1+\frac{4qD}%
{p^{2}v^{2}}s}}\right)   &  =\frac{p^{2}v^{2}}{4qD}\frac{1}{2\pi i}%
\int_{\left(  1+\frac{4qD}{p^{2}v^{2}}c\right)  -i\infty}^{\left(
1+\frac{4qD}{p^{2}v^{2}}c\right)  +i\infty}e^{-\xi\sqrt{\gamma}}e^{\frac
{p^{2}v^{2}}{4qD}\left(  \gamma-1\right)  t^{\prime}}d\gamma\nonumber\\
&  =\frac{p^{2}v^{2}}{4qD}e^{-\theta}\frac{1}{2\pi i}\int_{g-i\infty
}^{g+i\infty}e^{-\xi\sqrt{\gamma}}e^{\theta\gamma}d\gamma,\label{3.22a}\\
\text{where \ \ \ \ \ \ \ \ \ \ \ }g  &  =\left(  1+\frac{4qD}{p^{2}v^{2}%
}c\right)  \text{ \ and }\theta=\frac{p^{2}v^{2}}{4qD}t^{\prime}.
\label{3.22b}%
\end{align}
\medskip Rewriting the integral in function notation, we obtain\medskip%
\begin{equation}
\mathcal{L}^{-1}\left(  e^{-\frac{pvx^{\prime}}{2qD}\sqrt{1+\frac{4qD}%
{p^{2}v^{2}}s}}\right)  =\frac{p^{2}v^{2}}{4qD}e^{-\theta}\mathcal{L}%
^{-1}\left(  e^{-\xi\sqrt{\gamma}}\right)  . \label{3.23}%
\end{equation}
\medskip From the tables\cite{CRC1996} (p. 563, \#74) the corresponding
inverse Laplace transform is\medskip%
\begin{align}
\mathcal{L}^{-1}\left(  e^{-\xi\sqrt{\gamma}}\right)   &  =\frac{\xi}%
{2\theta\sqrt{\pi\theta}}e^{-\frac{\xi^{2}}{4\theta}},\label{3.24a}\\
\mathcal{L}^{-1}\left(  e^{-\frac{pvx^{\prime}}{2qD}\sqrt{1+\frac{4qD}%
{p^{2}v^{2}}s}}\right)   &  =\left(  \frac{p^{2}v^{2}}{4qD}e^{-\theta}\right)
\left(  \frac{\xi}{2\theta\sqrt{\pi\theta}}e^{-\frac{\xi^{2}}{4\theta}%
}\right)  . \label{3.24b}%
\end{align}
\medskip Substituting this result in (\ref{3.19}) and replacing $e^{-\frac
{pvx^{\prime}}{2qD}}$ with $e^{-\xi}$, we obtain\medskip\ $\ $%
\begin{equation}
\Psi\left(  x^{\prime},t^{\prime}\right)  =Q\left(  t^{\prime}\right)  \ast
e^{-\xi}\left(  \frac{p^{2}v^{2}}{4qD}e^{-\theta}\right)  \left(  \frac{\xi
}{2\theta\sqrt{\pi\theta}}e^{-\frac{\xi^{2}}{4\theta}}\right)  . \label{3.25}%
\end{equation}
\medskip Simplifying and then returning to our unscaled space and time
variables, we get\medskip\
\begin{align}
\Psi\left(  x^{\prime},t^{\prime}\right)   &  =Q\left(  t^{\prime}\right)
\ast\frac{p^{2}v^{2}}{2qD}\frac{\xi}{4\theta\sqrt{\pi\theta}}e^{-\frac{\left(
2\theta+\xi\right)  ^{2}}{4\theta}}\label{3.26a}\\
&  =Q\left(  t^{\prime}\right)  \ast\frac{x^{\prime}}{t^{\prime}\sqrt{4\pi
qDt^{\prime}}}e^{-\frac{\left(  x^{\prime}+pvt^{\prime}\right)  ^{2}%
}{4qDt^{\prime}}}. \label{3.26b}%
\end{align}
\medskip Finally, writing the convolution as an integral, we obtain\medskip%
\begin{equation}
\Psi\left(  x^{\prime},t^{\prime}\right)  =\int_{0}^{t^{\prime}}%
\frac{x^{\prime}}{\tau\sqrt{4\pi qD\tau}}e^{-\frac{\left(  x^{\prime}%
+pv\tau\right)  ^{2}}{4qD\tau}}Q\left(  t^{\prime}-\tau\right)  d\tau.
\label{3.27}%
\end{equation}
\smallskip(Here $\tau$ is just a dummy variable for the integration.)

By the mean value theorem there is some value of time $t_{\ast}^{\prime},$
lying in the interval $\left(  0,t^{\prime}\right)  ,$ such that$\medskip$%

\begin{equation}
\Psi\left(  x^{\prime},t^{\prime}\right)  =Q\left(  t_{\ast}^{\prime}\right)
\int_{0}^{t^{\prime}}\frac{x^{\prime}}{\tau\sqrt{4\pi qD\tau}}e^{-\frac
{\left(  x^{\prime}+pv\tau\right)  ^{2}}{4qD\tau}}d\tau. \label{3.28}%
\end{equation}
\medskip Using definitions (\ref{3.20a}) and (\ref{3.22b}), we switch back to
the unitless variables to make the integration simpler $\left(  \frac{2qD}%
{pv}\xi=x^{\prime}\right.  $ , $\frac{4qD}{p^{2}v^{2}}\theta=t^{\prime}$ , and
$\left.  \frac{4qD}{p^{2}v^{2}}u=\tau\right)  $:\medskip\
\begin{align*}
\int_{0}^{t^{\prime}}\frac{x^{\prime}}{\tau\sqrt{4\pi qD\tau}}e^{-\frac
{\left(  x^{\prime}+pv\tau\right)  ^{2}}{4qD\tau}}d\tau &  =\frac{\xi}%
{2\sqrt{\pi}}\int_{0}^{\theta}\frac{1}{\sqrt{u^{3}}}e^{-\frac{\left(
2u+\xi\right)  ^{2}}{4u}}du\\
&  =\frac{1}{2}\left[  e^{-2\xi}\left(  1+\operatorname{erf}\left(
\frac{2\theta-\xi}{2\sqrt{\theta}}\right)  \right)  +\operatorname{erfc}%
\left(  \frac{2\theta+\xi}{2\sqrt{\theta}}\right)  \right]  ,
\end{align*}
\medskip%
\begin{equation}
\Psi\left(  x^{\prime},t^{\prime}\right)  =Q\left(  t_{\ast}^{\prime}\right)
\frac{1}{2}\left[  e^{-\frac{pvx^{\prime}}{qD}}\left(  1+\operatorname{erf}%
\left(  \frac{pvt^{\prime}-x^{\prime}}{\sqrt{4qt^{\prime}D}}\right)  \right)
+\operatorname{erfc}\left(  \frac{pvt^{\prime}+x^{\prime}}{\sqrt{4qDt^{\prime
}}}\right)  \right]  . \label{3.29}%
\end{equation}
\medskip However, from the way we defined the problem, the total integral of
$\ \Psi\left(  x^{\prime},t^{\prime}\right)  $ is the total probability of an
uninhibited flow passing $x=a$ over all $x$ and is to equal $Q\left(
t^{\prime}\right)  .$Integrating (\ref{3.28}) over all $x^{\prime}$, we
have\medskip\
\begin{equation}
\int_{0}^{\infty}\Psi\left(  x^{\prime},t^{\prime}\right)  dx^{\prime
}=Q\left(  t^{\prime}\right)  . \label{3.30}%
\end{equation}
\medskip Therefore,%
\begin{equation}
Q\left(  t^{\prime}\right)  =Q\left(  t_{\ast}^{\prime}\right)  \int
_{0}^{\infty}\int_{0}^{t^{\prime}}\frac{x^{\prime}}{\tau\sqrt{4\pi qD\tau}%
}e^{-\frac{\left(  x^{\prime}+pv\tau\right)  ^{2}}{4qD\tau}}d\tau dx^{\prime}.
\label{3.31}%
\end{equation}
\medskip If we take the integral in (\ref{3.31}) as the probability
distribution density for $Q\left(  t^{\prime}\right)  $. then we need a
normalizing factor $\left(  h\left(  t^{\prime}\right)  \right)  $ to make it
integrate to unity.\medskip\
\begin{equation}
1=h\left(  t^{\prime}\right)  \int_{0}^{\infty}\int_{0}^{t^{\prime}}%
\frac{x^{\prime}}{\tau\sqrt{4\pi qD\tau}}e^{-\frac{\left(  x^{\prime}%
+pv\tau\right)  ^{2}}{4qD\tau}}d\tau dx^{\prime} \label{3.32}%
\end{equation}
\medskip Again we switch to the non-dimensional variables before integrating
$\left(  \frac{2qD}{pv}\xi=x^{\prime}\right.  $ , $\frac{4qD}{p^{2}v^{2}%
}\theta=t^{\prime}$ , and $\left.  \frac{4qD}{p^{2}v^{2}}u=\tau\right)  $, and
procced as follows\footnote{Note that $\operatorname{erf}\left(  x\right)
=\frac{2}{\pi}\int_{0}^{x}e^{-t^{2}}dt$ and $\operatorname{erfc}\left(
x\right)  =\frac{2}{\pi}\int_{x}^{\infty}e^{-t^{2}}dt=1-\operatorname{erf}%
\left(  x\right)  .$}\medskip%
\begin{align*}
1  &  =h\left(  t^{\prime}\right)  \frac{qD}{\sqrt{\pi}pv}\int_{0}^{\infty
}\int_{0}^{\theta}\xi u^{-\frac{3}{2}}e^{-\frac{1}{4u}\left(  2u+\xi\right)
^{2}}dud\xi\\
&  =h\left(  t^{\prime}\right)  \frac{qD}{\sqrt{\pi}pv}\int_{0}^{\theta
}u^{-\frac{3}{2}}\left(  \int_{0}^{\infty}\xi e^{-\frac{1}{4u}\left(
2u+\xi\right)  ^{2}}d\xi\right)  du\\
&  =h\left(  t^{\prime}\right)  \frac{2qD}{pv}\int_{0}^{\theta}\left(
\frac{1}{\sqrt{\pi u}}e^{-u}-\operatorname{erfc}\left(  \sqrt{u}\right)
\right)  du\\
&  =h\left(  t^{\prime}\right)  \frac{2qD}{pv}\left(  \frac{1}{2}%
\operatorname{erf}\left(  \sqrt{\theta}\right)  -\theta\operatorname{erfc}%
\left(  \sqrt{\theta}\right)  +\sqrt{\frac{\theta}{\pi}}e^{-\theta}\right)  .
\end{align*}%
\[
1=h\left(  t^{\prime}\right)  \left[  \frac{qD}{pv}\operatorname{erf}\left(
\sqrt{\frac{p^{2}v^{2}t^{\prime}}{4qD}}\right)  -\frac{pvt^{\prime}}%
{2}\operatorname{erfc}\left(  \sqrt{\frac{p^{2}v^{2}t^{\prime}}{4qD}}\right)
+\sqrt{\frac{qDt^{\prime}}{\pi}}e^{-\frac{p^{2}v^{2}}{4qD}t^{\prime}}\right]
.
\]
\medskip Solving for $h\left(  t^{\prime}\right)  ,$ we have\medskip\
\begin{equation}
h\left(  t^{\prime}\right)  =\left(  \frac{qD}{pv}\operatorname{erf}\left(
\frac{pvt^{\prime}}{\sqrt{4qDt^{\prime}}}\right)  -\frac{pvt^{\prime}}%
{2}\operatorname{erfc}\left(  \frac{pvt^{\prime}}{\sqrt{4qDt^{\prime}}%
}\right)  +\sqrt{\frac{qDt^{\prime}}{\pi}}e^{-\frac{p^{2}v^{2}t^{\prime}}%
{4qD}}\right)  ^{-1}. \label{3.33}%
\end{equation}
\medskip Note that from (\ref{3.31}) and (\ref{3.32}), we have\medskip\
\begin{equation}
Q\left(  t_{\ast}^{\prime}\right)  =h\left(  t^{\prime}\right)  Q\left(
t^{\prime}\right)  . \label{3.34}%
\end{equation}
\medskip Putting this together with (\ref{3.29}) and (\ref{3.33}), we
get\medskip%
\begin{equation}
\Psi\left(  x^{\prime},t^{\prime}\right)  =Q\left(  t^{\prime}\right)  \left(
\frac{e^{-\frac{pvx^{\prime}}{qD}}\left(  1+\operatorname{erf}\left(
\frac{pvt^{\prime}-x^{\prime}}{\sqrt{4qDt^{\prime}}}\right)  \right)
+\operatorname{erfc}\left(  \frac{pvt^{\prime}+x^{\prime}}{\sqrt{4qDt^{\prime
}}}\right)  }{\frac{2qD}{pv}\operatorname{erf}\left(  \frac{pvt^{\prime}%
}{\sqrt{4qDt^{\prime}}}\right)  -pvt^{\prime}\operatorname{erfc}\left(
\frac{pvt^{\prime}}{\sqrt{4qDt^{\prime}}}\right)  +\sqrt{\frac{4qDt^{\prime}%
}{\pi}}e^{-\frac{p^{2}v^{2}t^{\prime}}{4qD}}}\right)  \label{3.35}%
\end{equation}
\medskip\ As we did with the distribution parallel to the $y$-axis, we will
need to add reflections of this distribution to get the correct image within
the boundaries of the enclosure. \ The presence of $x^{\prime}$ terms in the
argument of the error functions will make for very messy expressions. To
simplify (a little bit), we can approximate this expression as just an
exponential. We obtain the approximation by evaluating the above at
$x^{\prime}=0$ and using the result for the parametrization of the
exponential. As it turns out this is just the normalizing factor, $h,$ that we
obtained above:\medskip%
\begin{subequations}
\label{3.36}%
\begin{align}
\Psi\left(  x^{\prime},t^{\prime}\right)   &  \approx Q\left(  t^{\prime
}\right)  h\left(  t^{\prime}\right)  e^{-h\left(  t^{\prime}\right)
x^{\prime}},\label{3.36a}\\
\text{where \ }h\left(  t^{\prime}\right)   &  =\left(  \frac{qD}%
{pv}\operatorname{erf}\left(  \frac{pvt^{\prime}}{\sqrt{4qDt^{\prime}}%
}\right)  -\frac{pvt^{\prime}}{2}\operatorname{erfc}\left(  \frac{pvt^{\prime
}}{\sqrt{4qDt^{\prime}}}\right)  +\sqrt{\frac{qDt^{\prime}}{\pi}}%
e^{-\frac{p^{2}v^{2}t^{\prime}}{4qD}}\right)  ^{-1}. \label{3.36b}%
\end{align}
\medskip Finally, returning to the original coordinate system, we
have\medskip\
\end{subequations}
\begin{equation}
\Psi\left(  x,t\right)  =Q\left(  t\right)  \left(  \frac{e^{-\frac{pv\left(
a-x\right)  }{qD}}\left(  1+\operatorname{erf}\left(  \frac{pvt-\left(
a-x\right)  }{\sqrt{4qDt}}\right)  \right)  +\operatorname{erfc}\left(
\frac{pvt+\left(  a-x\right)  }{\sqrt{4qDt}}\right)  }{\frac{2qD}%
{pv}\operatorname{erf}\left(  \frac{pvt}{\sqrt{4qDt}}\right)
-pvt\operatorname{erfc}\left(  \frac{pvt}{\sqrt{4qDt}}\right)  +\sqrt
{\frac{4qDt}{\pi}}e^{-\frac{p^{2}v^{2}t}{4qD}}}\right)  . \label{3.37}%
\end{equation}
\medskip\ The approximation we use is\medskip\
\begin{equation}
\Psi\left(  x,t\right)  \approx Q\left(  t\right)  h\left(  t\right)
e^{h\left(  t\right)  \left(  x-a\right)  }. \label{3.38}%
\end{equation}
\medskip

\begin{figure}[th]
\centering\includegraphics[width=6in,height=2.2in]{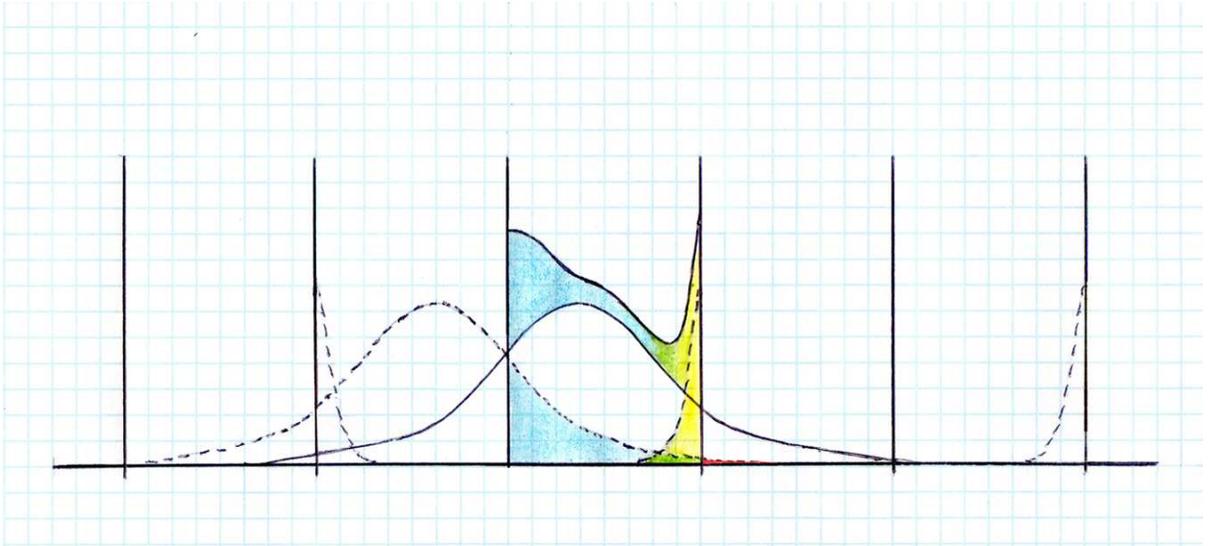}%
\caption{The primary distribution crosses the boundary at $x=0$ and $x=a$. The
left tail of the pdf can be replaced within the boundaries of the problem by
adding the tail of its "reflected" distribution centered at $x=-x_{0}-vt$ on
the left. To the right of $a$ we integrate the area of the tails of the
original distribution \textit{and} its image. This area (colored yellow) is
distributed within the boundaries as an exponential pdf with its maximum value
at $x=a$ . This distribution is fixed in location and has infinite images as
in the case for the $y$-axis. Each of these is superposed, as before, within
the practical boundary and is shown as blended colors over the primary
distribution.}%
\label{Figure:5}%
\end{figure}

To find the area redistributed, we scale $\Psi$ by the total area of the
original PDE \ $X_{0}$\ and its image $X_{-1},$ that is to the right of the
right most boundary at $x=a.$ The image $\Psi_{-1}$ is the only image, since
we do not treat the right boundary as reflecting the original (traveling)
distribution (see \ref{Figure:5}). The original distribution and its image is
given by\medskip\
\begin{equation}%
\begin{array}
[c]{l}%
X_{0}\left(  x,t\right)  =\frac{1}{\sqrt{4\pi qDt}}e^{-\frac{\left(  x-\left(
x_{0}+pvt\right)  \right)  ^{2}}{4qDt}},\\
X_{-1}\left(  x,t\right)  =\frac{1}{\sqrt{4\pi qDt}}e^{-\frac{\left(
x+\left(  x_{0}+pvt\right)  \right)  ^{2}}{4qDt}}.
\end{array}
\label{3.39}%
\end{equation}
\medskip The total area redistributed is then\medskip%
\begin{align}
Q\left(  t\right)   &  =\int_{a}^{\infty}\left(  X_{0}\left(  x,t\right)
+X_{-1}\left(  x,t\right)  \right)  dx\nonumber\\
&  =\frac{1}{\sqrt{4\pi qDt}}\int_{a}^{\infty}\left(  e^{-\frac{\left(
x-\left(  x_{0}+pvt\right)  \right)  ^{2}}{4qDt}}+e^{-\frac{\left(  x+\left(
x_{0}+pvt\right)  \right)  ^{2}}{4qDt}}\right)  dx\nonumber\\
&  =1+\frac{1}{2}\left(  \operatorname{erf}\left(  \frac{x_{0}-a+pvt}%
{\sqrt{4qDt}}\right)  -\operatorname{erf}\left(  \frac{\left(  x_{0}%
+a+pvt\right)  }{\sqrt{4qDt}}\right)  \right)  . \label{3.40}%
\end{align}
$\allowbreak$ The scaled distribution is\medskip%
\begin{equation}
\Psi\left(  x,t\right)  \approx Q\left(  t\right)  h\left(  t\right)
e^{-h\left(  t\right)  \left(  a-x\right)  }. \label{3.41}%
\end{equation}
\medskip The next part of the problem is\ that the (left) tail of this
distribution goes past the $x=0$ boundary. \ However, the $\Psi$ distribution,
though growing, does not move. \ We can then treat this distribution as we did
the distribution parallel to the $y$-axis. As in that case, we take the part
of each reflected image that lies between $x=0$ and $x=a$ and add it to the
original distribution. \ If the $\Psi$\ distribution is anchored at $a,$ then
the distribution and images are anchored at $\left(  2m+1\right)  a$ for
$m=\left\{  0,\pm1,\pm2,\pm3,\ldots\right\}  $.

\medskip\medskip%

\begin{align}
\sum_{m=0}^{\infty}\Psi_{m}\left(  x,t\right)   &  =\sum_{m=0}^{\infty}\left(
Q\left(  t\right)  h\left(  t\right)  e^{h\left(  t\right)  \left(  x-\left(
2m+1\right)  a\right)  }+Q\left(  t\right)  h\left(  t\right)  e^{h\left(
t\right)  \left(  -x-\left(  2m+1\right)  a\right)  }\right)  \text{ \ for
}\ 0\leq x\leq a\nonumber\\
&  =Q\left(  t\right)  h\left(  t\right)  \sum_{m=0}^{\infty}\left(
e^{h\left(  t\right)  x}e^{-h\left(  t\right)  a}e^{-2mh\left(  t\right)
a}+e^{-h\left(  t\right)  x}e^{-h\left(  t\right)  a}e^{-2mh\left(  t\right)
a}\right) \nonumber\\
&  =Q\left(  t\right)  h\left(  t\right)  e^{-h\left(  t\right)  a}\left(
e^{h\left(  t\right)  x}+e^{-h\left(  t\right)  x}\right)  \sum_{m=0}^{\infty
}e^{-2mh\left(  t\right)  a}v \label{3.42}%
\end{align}
\medskip Note that $\left(  e^{h\left(  t\right)  x}+e^{-h\left(  t\right)
x}\right)  =2\cosh\left[  h\left(  t\right)  x\right]  ,$ so then\medskip\
\begin{equation}
\sum_{m=0}^{\infty}\Psi_{m}\left(  x,t\right)  =Q\left(  t\right)  2h\left(
t\right)  e^{-h\left(  t\right)  a}\cosh\left[  h\left(  t\right)  x\right]
\sum_{m=0}^{\infty}\left(  e^{-2h\left(  t\right)  a}\right)  ^{m}.
\label{3.43}%
\end{equation}
\medskip Note that $\sum_{m=0}^{\infty}\left(  e^{-2h\left(  t\right)
a}\right)  ^{m}$ is a geometric series, equal to $\ \left(  1-e^{-2h\left(
t\right)  a}\right)  ^{-1},$ hence,\medskip%
\begin{equation}
\sum_{m=0}^{\infty}\Psi_{m}\left(  x,t\right)  =Q\left(  t\right)  2h\left(
t\right)  \frac{e^{-h\left(  t\right)  a}}{1-e^{-2h\left(  t\right)  a}}%
\cosh\left[  h\left(  t\right)  x\right]  =Q\left(  t\right)  h\left(
t\right)  \frac{\cosh\left[  h\left(  t\right)  x\right]  }{\sinh\left[
h\left(  t\right)  a\right]  }, \label{3.44}%
\end{equation}
\medskip where the $\sinh$ term came from the fact that
\[
\frac{e^{-h\left(  t\right)  a}}{1-e^{-2h\left(  t\right)  a}}=\frac
{1}{e^{h\left(  t\right)  a}-e^{-h\left(  t\right)  a}}=\frac{1}{2\sinh\left(
h\left(  t\right)  a\right)  }.
\]
\medskip

Putting the pieces together, we get\medskip\
\begin{align}
&  \left(  X_{0}\left(  x,t\right)  +X_{-1}\left(  x,t\right)  \right)
+\sum_{m=0}^{\infty}\Psi_{m}\left(  x,t\right) \nonumber\\
&  =\frac{1}{\sqrt{4\pi qDt}}\left(  e^{-\frac{\left(  x-x_{0}-pvt\right)
^{2}}{4qDt}}+e^{-\frac{\left(  x+x_{0}+pvt\right)  ^{2}}{4qDt}}\right)
+Q\left(  t\right)  h\left(  t\right)  \frac{\cosh\left(  h\left(  t\right)
x\right)  }{\sinh\left(  h\left(  t\right)  a\right)  }. \label{3.45}%
\end{align}
\medskip Finally, multiplying the $x$ and $y$ oriented distributions, we have
the following solution on \ \ $0<x<a$\medskip%
\begin{subequations}
\begin{align}
P\left(  x,y,t\right)   &  =\left(  \frac{1}{\sqrt{4\pi qDt}}\left(
e^{-\frac{\left(  x-x_{0}-pvt\right)  ^{2}}{4qDt}}+e^{-\frac{\left(
x+x_{0}+pvt\right)  ^{2}}{4qDt}}\right)  +Q\left(  t\right)  h\left(
t\right)  \frac{\cosh\left(  h\left(  t\right)  x\right)  }{\sinh\left(
h\left(  t\right)  a\right)  }\right) \nonumber\\
&  \times\frac{1}{\sqrt{4\pi qDt}}\sum_{k=-\infty}^{\infty}\left(
e^{-\frac{\left(  y-y_{0}+2kb\right)  ^{2}}{4qDt}}+e^{-\frac{\left(
y+y_{0}+\left(  2k+1\right)  b\right)  ^{2}}{4qDt}}\right)  ,\label{3.46a}\\
\text{where }Q\left(  t\right)   &  =1+\frac{1}{2}\left(  \operatorname{erf}%
\left(  \frac{x_{0}-a+pvt}{\sqrt{4qDt}}\right)  -\operatorname{erf}\left(
\frac{\left(  x_{0}+a+pvt\right)  }{\sqrt{4qDt}}\right)  \right)
\label{3.46b}\\
\text{and \ }h\left(  t\right)   &  =\left(  \frac{qD}{pv}\operatorname{erf}%
\left(  \frac{pvt}{\sqrt{4qDt}}\right)  -\frac{pvt}{2}\operatorname{erfc}%
\left(  \frac{pvt}{\sqrt{4qDt}}\right)  +\sqrt{\frac{qDt}{\pi}}e^{-\frac
{p^{2}v^{2}t}{4qD}}\right)  ^{-1}. \label{3.46c}%
\end{align}
$\medskip$

As $t$ gets very large, we have the following limits

$\medskip$%
\end{subequations}
\begin{equation}%
\begin{array}
[c]{l}%
\lim_{t\rightarrow\infty}\frac{1}{\sqrt{4\pi qDt}}\left(  e^{-\frac{\left(
x-x_{0}-pvt\right)  ^{2}}{4qDt}}+e^{-\frac{\left(  x+x_{0}+pvt\right)  ^{2}%
}{4qDt}}\right)  =0.\\
\lim_{t\rightarrow\infty}\left(  \frac{1}{2}\left(  2+\operatorname{erf}%
\left(  \frac{pvt-a+x_{0}}{\sqrt{4qDt}}\right)  -\operatorname{erf}\left(
\frac{pvt+a+x_{0}}{\sqrt{4qDt}}\right)  \right)  \right)  =1.\\
\lim_{t\rightarrow\infty}\frac{1}{\sqrt{4\pi qDt}}\sum_{k=-\infty}^{\infty
}\left(  e^{-\frac{\left(  y-y_{0}+2kb\right)  ^{2}}{4qDt}}+e^{-\frac{\left(
y+y_{0}+\left(  2k+1\right)  b\right)  ^{2}}{4qDt}}\right)  =\frac{1}{b}.\\
\lim_{t\rightarrow\infty}Q\left(  t\right)  =1.\\
\lim_{t\rightarrow\infty}h\left(  t\right)  =\frac{pv}{qD}.\\
\lim_{t\rightarrow\infty}\frac{\cosh\left(  h\left(  t\right)  x\right)
}{\sinh\left(  h\left(  t\right)  a\right)  }=\frac{\cosh\left(  \frac
{pvx}{qD}\right)  }{\sinh\left(  \frac{pva}{qD}\right)  }.
\end{array}
\label{3.47}%
\end{equation}
\medskip%
\begin{equation}
P\left(  x,y,\infty\right)  =\frac{pv}{qDb}\frac{\cosh\left(  \frac{pvx}%
{qD}\right)  }{\sinh\left(  \frac{pva}{qD}\right)  }\text{ \ \ \ \ \ }%
\left\{
\begin{array}
[c]{c}%
0<x<a\\
\frac{-b}{2}<y<\frac{-b}{2}%
\end{array}
\right\}  . \label{3.48}%
\end{equation}
\medskip This is the steady-state solution for\ the probability distribution.

\subsection{Median time to first arrival at goal.}

We can also use the expression for $Q\left(  t\right)  $ to obtain the median
time to reach the goal. This is the time, $t_{M},$defined when $Q\left(
t_{M}\right)  =1/2.$%
\begin{align}
Q\left(  t_{M}\right)   &  =\frac{1}{2}=1+\frac{1}{2}\left(
\operatorname{erf}\left(  \frac{x_{0}-a+pvt_{M}}{\sqrt{4qDt_{M}}}\right)
-\operatorname{erf}\left(  \frac{\left(  x_{0}+a+pvt_{M}\right)  }%
{\sqrt{4qDt_{M}}}\right)  \right) \nonumber\\
1  &  =\operatorname{erf}\left(  \frac{\left(  x_{0}+a+pvt_{M}\right)  }%
{\sqrt{4qDt_{M}}}\right)  -\operatorname{erf}\left(  \frac{x_{0}-a+pvt_{M}%
}{\sqrt{4qDt_{M}}}\right)  \label{3.48a}%
\end{align}
This, of course, has to be solved numerically.\medskip

\subsection{Summary\medskip}

Summarizing our results, the probability distribution with the original
parametrization is\medskip\
\begin{align}
P\left(  x,y,t\right)  dxdy  &  =\left(  \frac{1}{\sqrt{4\pi qDt}}\left(
e^{-\frac{\left(  x-x_{0}-pvt\right)  ^{2}}{4qDt}}+e^{-\frac{\left(
x+x_{0}+pvt\right)  ^{2}}{4qDt}}\right)  +Q\left(  t\right)  h\left(
t\right)  \frac{\cosh\left(  h\left(  t\right)  x\right)  }{\sinh\left(
h\left(  t\right)  a\right)  }\right) \nonumber\\
&  \times\frac{1}{\sqrt{4\pi qDt}}\sum_{k=-\infty}^{\infty}\left(
e^{-\frac{\left(  y-y_{0}+2kb\right)  ^{2}}{4qDt}}+e^{-\frac{\left(
y+y_{0}+\left(  2k+1\right)  b\right)  ^{2}}{4qDt}}\right)  dxdy\label{3.49}\\
\text{where \ }Q\left(  t\right)   &  =1+\frac{1}{2}\left(  \operatorname{erf}%
\left(  \frac{x_{0}-a+pvt}{\sqrt{4qDt}}\right)  -\operatorname{erf}\left(
\frac{\left(  x_{0}+a+pvt\right)  }{\sqrt{4qDt}}\right)  \right) \nonumber\\
\text{and \ \ \ \ }h\left(  t\right)   &  =\left(  \frac{qD}{pv}%
\operatorname{erf}\left(  \frac{pvt}{\sqrt{4qDt}}\right)  -\frac{pvt}%
{2}\operatorname{erfc}\left(  \frac{pvt}{\sqrt{4qDt}}\right)  +\sqrt
{\frac{qDt}{\pi}}e^{-\frac{p^{2}v^{2}t}{2qD}}\right)  ^{-1},\nonumber
\end{align}
or rescaling time and distance, using the following transformations%
\[%
\begin{array}
[c]{ccc}%
x=\frac{2qD}{pv}\xi & \ x_{0}=\frac{2qD}{pv}\xi_{0} & a=\frac{2qD}{pv}\alpha\\
y=\frac{2qD}{pv}\zeta & \ y_{0}=\frac{2qD}{pv}\zeta_{0} & b=\frac{2qD}%
{pv}\beta\\
& t=\frac{4qD}{p^{2}v^{2}}\theta. &
\end{array}
\]

we also have,%
\begin{align*}
P\left(  \xi,\zeta,\theta\right)  d\xi d\zeta &  =\left(  \frac{1}{2\sqrt
{\pi\theta}}\left(  e^{-\frac{\left(  \xi-\xi_{0}-2\theta\right)  ^{2}%
}{4\theta}}+e^{-\frac{\left(  \xi+\xi_{0}+2\theta\right)  ^{2}}{4\theta}%
}\right)  +Q\left(  \theta\right)  h\left(  \theta\right)  \frac{\cosh\left(
h\left(  \theta\right)  \xi\right)  }{\sinh\left(  h\left(  \theta\right)
\alpha\right)  }\right) \\
&  \times\frac{1}{2\sqrt{\pi\theta}}\sum_{k=-\infty}^{\infty}\left(
e^{-\frac{\left(  \zeta-\zeta_{0}+2k\beta\right)  ^{2}}{4\theta}}%
+e^{-\frac{\left(  \zeta+\zeta_{0}+\left(  2k+1\right)  \beta\right)  ^{2}%
}{4\theta}}\right)  d\xi d\zeta\\
\text{where \ }Q\left(  \theta\right)   &  =1+\frac{1}{2}\operatorname{erf}%
\left(  \frac{\xi_{0}-\alpha+2\theta}{2\sqrt{\theta}}\right)  -\frac{1}%
{2}\operatorname{erf}\left(  \frac{\xi_{0}+\alpha+2\theta}{2\sqrt{\theta}%
}\right) \\
\text{and \ \ \ \ }h\left(  \theta\right)   &  =2\left(  \operatorname{erf}%
\left(  \sqrt{\theta}\right)  -2\theta\operatorname{erfc}\left(  \sqrt{\theta
}\right)  +2\sqrt{\frac{\theta}{\pi}}e^{-\theta}\right)  ^{-1}.
\end{align*}
\medskip

\begin{figure}[th]
\centering\includegraphics[width=5.5in,height=1.87in]{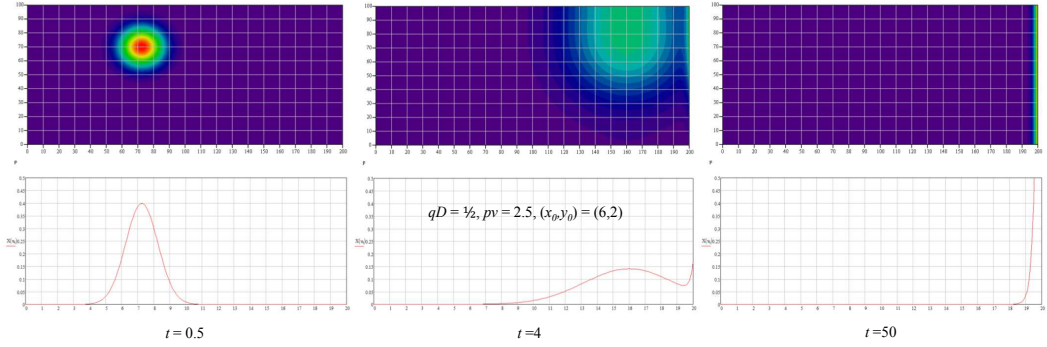}%
\caption{Evolution of the probability distribution at three points in time. }%
\label{Figure:6}%
\end{figure}

\section{Calculation of the probability of reaching the goal by a given point
in time, and the probability of placing at that time.\medskip}

The probability of fish species 1 reaching a point along the right-hand
boundary (located at $x=a$) by time $t$ is given by $Q_{1}\left(  t\right)  $
(the subscript indicating species 1)\medskip%
\begin{align}
&  Q_{1}\left(  t\right)  =\int_{a}^{\infty}\frac{1}{\sqrt{4\pi q_{1}D_{1}t}%
}\left(  e^{-\frac{\left(  x-x_{0}-p_{1}v_{1}t\right)  ^{2}}{4q_{1}D_{1}t}%
}+e^{-\frac{\left(  x+x_{0}+p_{1}v_{1}t\right)  ^{2}}{4q_{1}D_{1}t}}\right)
dx\nonumber\\
&  =1+\frac{1}{2}\left(  \operatorname{erf}\left(  \frac{p_{1}v_{1}t+x_{0}%
-a}{\sqrt{4q_{1}D_{1}t}}\right)  -\operatorname{erf}\left(  \frac{p_{1}%
v_{1}t+x_{0}+a}{\sqrt{4q_{1}D_{1}t}}\right)  \right)  . \label{4.01}%
\end{align}
\medskip The probability of\ fish species 2 (subscript 2) not reaching a point
along the right hand boundary (located at $x=a$) by time $t$ is given by
$1-Q_{2}\left(  t\right)  $\medskip%
\begin{equation}
1-Q_{2}\left(  t\right)  =\frac{1}{2}\left(  \operatorname{erf}\left(
\frac{p_{2}v_{2}t+x_{0}+a}{\sqrt{4q_{2}D_{2}t}}\right)  -\operatorname{erf}%
\left(  \frac{p_{2}v_{2}t+x_{0}-a}{\sqrt{4q_{2}D_{2}t}}\right)  \right)  .
\label{4.02}%
\end{equation}
\medskip Thus, the probability of fish species 1, beating fish species 2 by
time $t,$ is given by the product $\ Q_{1}\left(  t\right)  \left(
1-Q_{2}\left(  t\right)  \right)  :$\medskip%
\begin{align}
&  \left[  1+\frac{1}{2}\left(  \operatorname{erf}\left(  \frac{p_{1}%
v_{1}t+x_{0}-a}{\sqrt{4q_{1}D_{1}t}}\right)  -\operatorname{erf}\left(
\frac{p_{1}v_{1}t+x_{0}+a}{\sqrt{4q_{1}D_{1}t}}\right)  \right)  \right]
\nonumber\\
&  \times\frac{1}{2}\left(  \operatorname{erf}\left(  \frac{p_{2}v_{2}%
t+x_{0}+a}{\sqrt{4q_{2}D_{2}t}}\right)  -\operatorname{erf}\left(  \frac
{p_{2}v_{2}t+x_{0}-a}{\sqrt{4q_{2}D_{2}t}}\right)  \right)  . \label{4.03}%
\end{align}
\medskip\begin{figure}[th]
\centering\includegraphics[width=5.5in,height=2.5in]{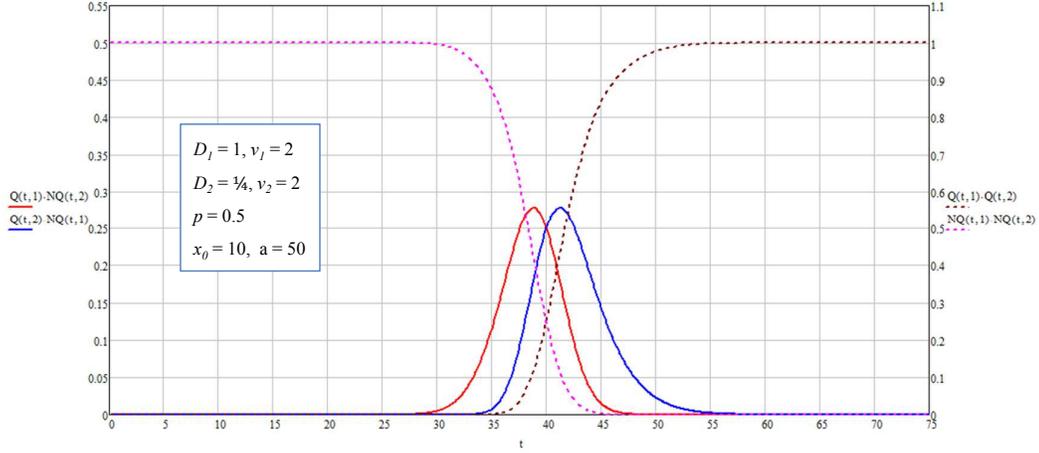}%
\caption{Competition between two species. Probability that a given species
will arrive first (left axis, solid curves). Probability that neither species
has arrived or that both have arrived (pink and brown dashed lines
respectively, scale on right axis). }%
\label{Figure:7}%
\end{figure}In general,he probability of fish species $k$ arriving first of
$n$ fish species is given by\medskip\
\begin{equation}
\ \frac{Q_{k}\left(  t\right)  }{\left(  1-Q_{k}\left(  t\right)  \right)
}\prod_{i=1}^{n}\left(  1-Q_{i}\left(  t\right)  \right)  . \label{4.04}%
\end{equation}
\medskip\begin{figure}[thth]
\centering\includegraphics[width=5.5in,height=2.75in]{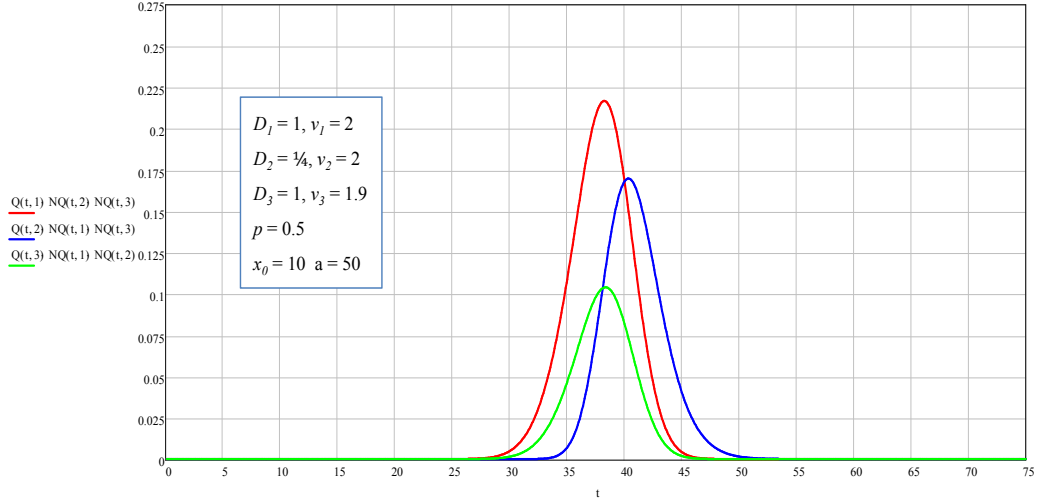}%
\caption{Competition between three species. Probability that a given species
will arrive first.}%
\label{Figure:8}%
\end{figure}The relative proportion of arrivals of species $k$ by time $t$ out
of all $n$ species is\medskip%
\begin{equation}
\ F\left(  t,n,k\right)  =\frac{N_{k}Q_{k}\left(  t\right)  }{\sum_{i=1}%
^{n}N_{i}Q_{i}\left(  t\right)  }, \label{4.05}%
\end{equation}
\medskip where $N_{i}$ is the total number of individuals in species
$i.$\medskip\begin{figure}[ththth]
\centering\includegraphics[width=5.5in,height=2.78in]{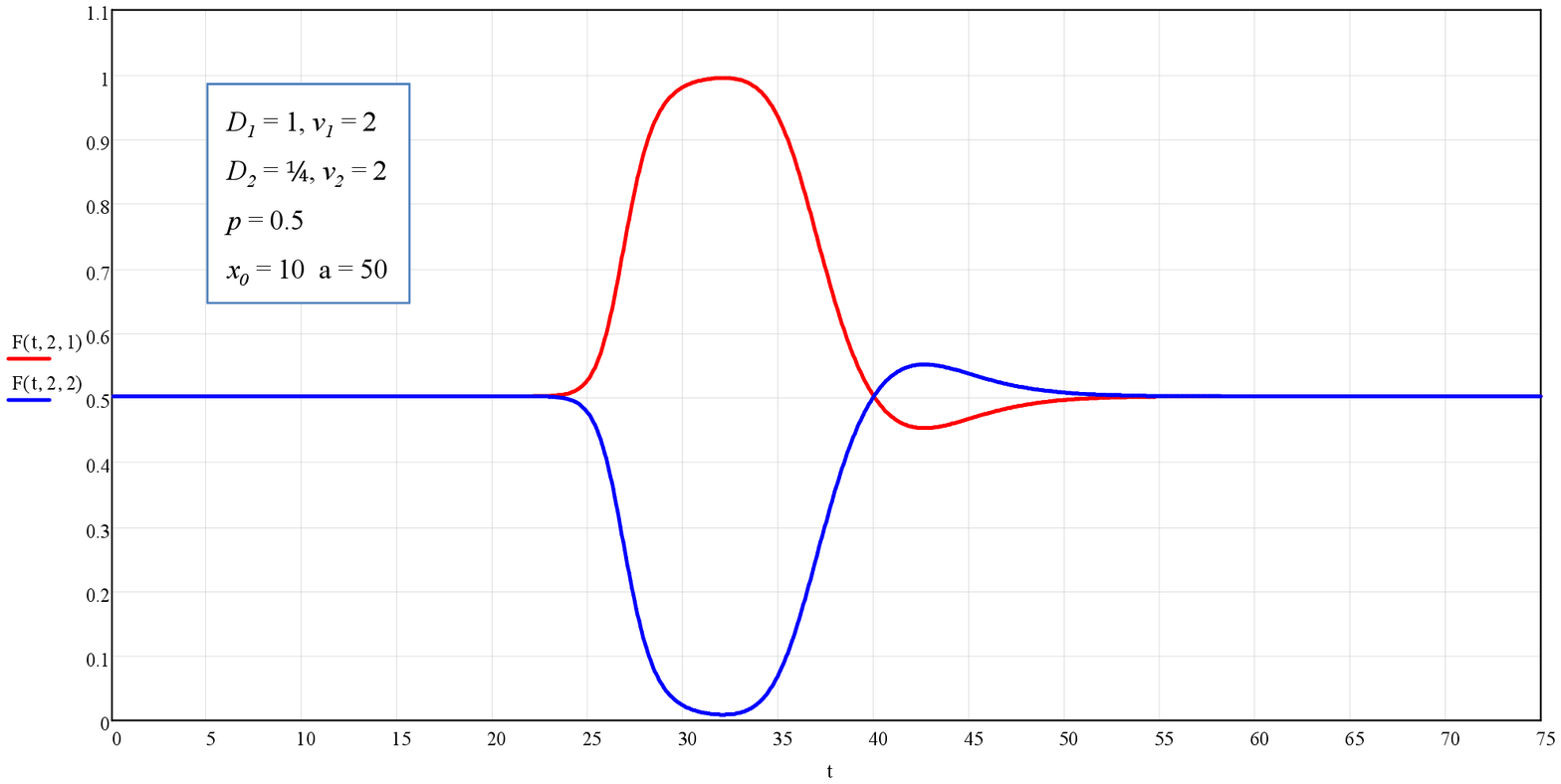}%
\caption{Competition between two species. Fraction of the total number of
individuals that members of a given species comprise. (If the probability of a
individual arriving is less than $10^{-4}$ then probability is set to
$10^{-4},$ this prevents division by very small numbers.)}%
\label{Figure:9}%
\end{figure}

\begin{figure}[th]
\centering\includegraphics[width=5.5in,height=2.78in]{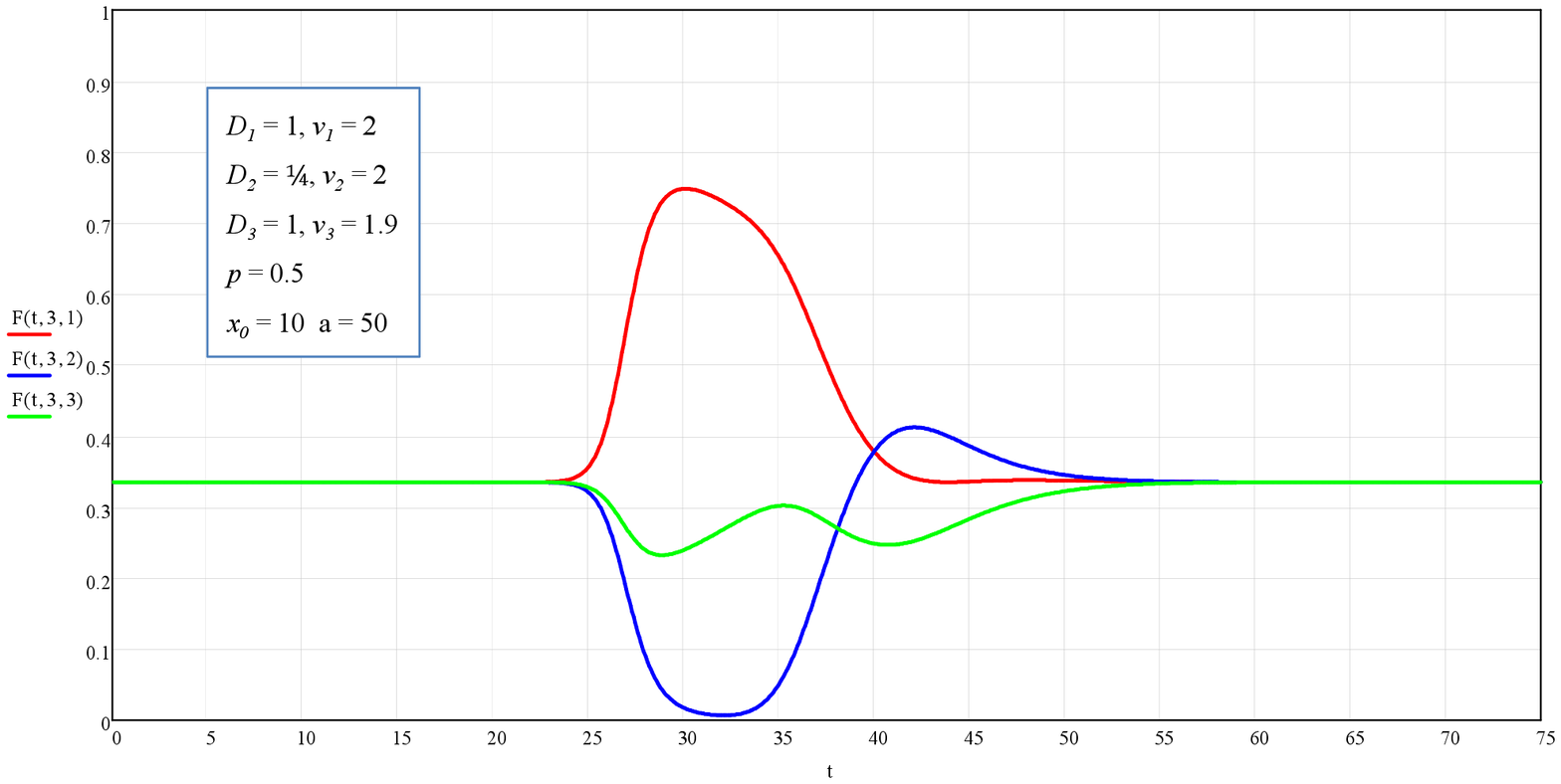}%
\caption{Competition between three species. Fraction of the total number of
individuals that members of a given species comprise. (If the probability of a
individual arriving is less than $10^{-4}$ then probability is set to
$10^{-4},$ this prevents division by very small numbers.)}%
\label{Figure:10}%
\end{figure}

We need to be a bit more specific with function arguments so that we can
specify order of arrival of distributions for placing second, third, and so
on. For the marginal distribution on the $x$ axis, we let\medskip\
\begin{equation}
P_{i}(x,t)=X_{0}\left(  x,t\right)  +X_{-1}\left(  x,t\right)  ,\text{\ \ }
\label{4.06}%
\end{equation}
\medskip and the probability of passing point $x$ is given by%
\begin{align}
Q_{i}\left(  x,t\right)   &  =\int_{x}^{\infty}P_{i}(\xi,t)d\xi\nonumber\\
&  =1+\frac{1}{2}\left(  \operatorname{erf}\left(  \frac{x_{0}-x+p_{i}v_{i}%
t}{\sqrt{4q_{i}D_{i}t}}\right)  -\operatorname{erf}\left(  \frac{x_{0}%
+x+p_{i}v_{i}t}{\sqrt{4q_{i}D_{i}t}}\right)  \right)  . \label{4.07}%
\end{align}
\medskip Second Place is\medskip\
\begin{equation}
\
{\displaystyle\sum_{i\neq k}}
\frac{\int_{a}^{\infty}P_{k}(x,t)Q_{i}\left(  x,t\right)  dx}{\left(
1-Q_{i}\left(  a,t\right)  \right)  }\frac{\prod_{j=1}^{n}\left(
1-Q_{j}\left(  t\right)  \right)  }{Q_{k}\left(  a,t\right)  \left(
1-Q_{k}\left(  a,t\right)  \right)  }. \label{4.08}%
\end{equation}
\medskip Third place is\medskip\
\begin{equation}
\
{\displaystyle\sum_{i\neq k}}
{\displaystyle\sum_{j\neq k\ i\neq j}}
\frac{\int_{a}^{\infty}P_{k}(x,t)Q_{i}\left(  x,t\right)  Q_{j}\left(
x,t\right)  dx}{\left(  1-Q_{i}\left(  a,t\right)  \right)  \left(
1-Q_{j}\left(  a,t\right)  \right)  }\frac{\prod_{h=1}^{n}\left(
1-Q_{h}\left(  t\right)  \right)  }{Q_{k}\left(  a,t\right)  \left(
1-Q_{k}\left(  a,t\right)  \right)  }. \label{4.09}%
\end{equation}
\medskip And so on.

\section{Discussion}

Though we have used examples of an individual's movements to illustrate the
rationale behind the development of equations, in truth the distributions and
movements described are a picture of the actions of an infinite number of
individuals each starting from the same conditions and acting under the same
rules. Thus, though an individual fish may swim here or there and never appear
to have any goal in mind, it may be found that, over time, the whole school
moves toward one destination as if the school itself had a goal and self
determination. Yet, if we knew ahead of time the arrow of that mass movement
and we broke down the components of the movements of each individual along
that axis and its perpendicular, we would find that the whole of the action
would be just an infinitesimal excess of desire for a slight majority of
individuals of movement in that one direction. It is thus for creatures as it
is for molecules of air in a breeze, though each moves in any direction, on
average they all move toward one direction. When confined to an enclosure,
they eventually pile up at the boundary toward which they tend. This changes
their distribution in space from a dispersing bell shaped curve to a steady
state exponential shaped curve, as the forces of diffusion and drift balance
one another.
\clearpage\newpage\medskip
\bibliographystyle{unsrtnat}
\bibliography{pde}

\bigskip
\end{document}